\documentclass[12pt]{amsart}
\usepackage{amsthm, amsmath, amssymb, graphicx, amsfonts}
\usepackage{natbib}

\graphicspath{{Figures/}}

\newtheorem{thm}{Theorem}[section]
\newtheorem{lem}[thm]{Lemma}
\newtheorem{prop}[thm]{Proposition}
\newtheorem{cor}[thm]{Corollary}

\newtheorem{rem}[thm]{Remark}
\newtheorem{expl}[thm]{Example}

\newcommand{\NN}{{\mathbb N}}
\newcommand{\FF}{{\mathbb F}}

\newcommand{\RR}{{\mathbb R}}

\newcommand{\argmax}{\mathop{\rm{argmax}\mathstrut}\nolimits}
\newcommand{\hypo}{\mathop{\rm{hypo}\mathstrut}\nolimits}
\newcommand{\cl}{\mathop{\mathrm{cl}\mathstrut}\nolimits}
\newcommand{\conv}{\mathop{\mathrm{conv}\mathstrut}\nolimits}

\newcommand{\eps}{\epsilon}

\begin{document}

\renewcommand{\baselinestretch}{1.2}  
\markright{
}
\markboth{\hfill{\footnotesize\rm BALABDAOUI, JANKOWSKI, PAVLIDES, SEREGIN, WELLNER
}\hfill}
{\hfill {\footnotesize\rm GRENANDER ESTIMATOR AT ZERO} \hfill}
\renewcommand{\thefootnote}{}
$\ $\par
\fontsize{10.95}{14pt plus.8pt minus .6pt}\selectfont
\vspace{0.8pc}
\centerline{\large\bf ON THE GRENANDER ESTIMATOR AT ZERO}
\vspace{.4cm}
\centerline{Fadoua Balabdaoui$^1$, Hanna Jankowski$^2$, Marios Pavlides$^3$,}
\centerline{Arseni Seregin$^4$, and Jon Wellner$^4$}
\vspace{.4cm}
\centerline{\it $^1$Universit\'{e} Paris-Dauphine, $^2$York University, $^3$Frederick University Cyprus}
\centerline{\it and \ $^4$University of Washington}
\vspace{.55cm}
\fontsize{9}{11.5pt plus.8pt minus .6pt}\selectfont

\begin{quotation}
\noindent {\it Abstract:}
We establish limit theory for the Grenander estimator of a monotone density near zero. In particular we consider the situation when the true density $f_0$ is unbounded at zero, with different rates of growth to infinity.  In the course of our study we develop new switching relations by use of tools from convex analysis.  The theory is applied to a problem involving mixtures.

\vspace{9pt}
\par
\end{quotation}\par

\fontsize{11}{14pt plus.8pt minus .6pt}\selectfont

\section{Introduction and Main Results}
\label{sec:IntroMainResults}

Let $X_1, \ldots , X_n$ be a sample from a decreasing density $f_0$ on $(0,\infty)$,
and let $\widehat{f}_n$ denote the Grenander estimator (i.e. the maximum likelihood estimator)
of $f_0$.  Thus $\widehat{f}_n \equiv \widehat{f}_n^L$ is the \textit{ left derivative}
of the least concave majorant $\widehat{F}_n$ of
the empirical distribution function $\FF_n$; see e.g.
\citet{MR0086459, MR0093415}, \citet{MR822052}, and \citet[chapter 8]{MR891874}.

The Grenander estimator $\widehat{f}_n$ is a uniformly consistent estimator
of $f_0$ on sets bounded away from $0$ if $f_0$ is continuous:
\begin{eqnarray*}
\sup_{x \ge c} | \widehat{f}_n (x) - f_0 (x) | \rightarrow_{a.s.} 0
\end{eqnarray*}
for each $c>0$.  It is also known that $\widehat{f}_n$ is consistent with respect to the $L_1 $
($\| p-q \|_1 \equiv \int | p(x) - q(x) | dx$) and
Hellinger
($h^2(p,q) \equiv 2^{-1} \int \left [ \sqrt{p(x)} -
 \sqrt{q(x)} \right ]^2 dx$) metrics:  that is,
$$
\| \widehat{f}_n - f_0 \|_1 \rightarrow_{a.s.} 0 \qquad \mbox{and}
\qquad  h(\widehat{f}_n , f_0 ) \rightarrow_{a.s.} 0;
$$
see e.g.
\citet[Theorem 8.3, page 144]{MR891874} and \citet{MR1212164}.

However, it is also known that $\widehat{f}_n (0) \equiv \widehat{f}_n (0+) $
is an inconsistent estimator of $f_0(0) \equiv f_0 (0+) = \lim_{x \searrow 0} f_0 (x)$, even when
$f_0(0) < \infty$.  In fact,
\citet{MR1243398}
showed that
\begin{eqnarray}
\widehat{f}_n (0) \rightarrow_d f_0 (0) \sup_{t > 0} \frac{\NN(t)}{t} \stackrel{d}{=} f_0 (0) \frac{1}{U}
\label{WoodroofeSunThmBoundedCase}
\end{eqnarray}
as $n \rightarrow \infty$ where $\NN$ is a standard Poisson process on $[0,\infty)$
and $U \sim \ $Uniform$(0,1)$.
\citet{MR1243398} introduced penalized estimators $\widetilde{f}_n $ of $f_0$ which yield consistency at $0$:
$\widetilde{f}_n (0) \rightarrow_p  f_0 (0)$.
\citet{MR2283391} study estimation of $f_0 (0)$ based on the Grenander estimator $\widehat{f}_n$ evaluated at
points of the form $t= c n^{-\gamma}$.   Among other things, they show that
$\widehat{f}_n(n^{-1/3}) \rightarrow_p f_0(0)$  if $| f_0 ' (0+)| > 0$.

Our view in this paper is that the inconsistency of $\widehat{f}_n (0)$ as an estimator of $f_0 (0)$ exhibited in
(\ref{WoodroofeSunThmBoundedCase})
can be regarded
as a simple consequence of the fact
that the class of all monotone decreasing densities on $(0,\infty)$ includes many densities $f$ which
are unbounded at $0$, so that $f(0) = \infty$, and the Grenander estimator $\widehat{f}_n$ simply has
difficulty deciding which is true, even when $f_0 (0) < \infty$.
From this perspective we would like to have answers to the following three questions
under some reasonable hypotheses concerning the growth
of $f_0 (x)$ as $x \searrow 0$:
\begin{list}{}
        {\setlength{\topsep}{3pt}
        \setlength{\parskip}{0pt}
        \setlength{\partopsep}{0pt}
        \setlength{\parsep}{0pt}
        \setlength{\itemsep}{1pt}
        \setlength{\leftmargin}{30pt}}
\item[\textbf{Q1:}]
How fast does $\widehat{f}_n (0)$ diverge as
$n \rightarrow \infty$?
\item[\textbf{Q2:}]
Do the stochastic processes $\{ b_n \widehat{f}_n (a_n t): \ 0 \le t \le c \} $ converge
for some sequences
$a_n$, $b_n$, and $c>0$?
\item[\textbf{Q3:}]
What is the behavior of the relative error
$$
\sup_{0 \le x \le c_n} \bigg | \frac{\widehat{f}_n (x)}{f_0 (x)} - 1 \bigg | 
$$
for some constant $c_n$?
\end{list}
\medskip

It turns out that answers to questions {\bf Q1} - {\bf Q3}
are intimately related to the limiting behavior of the minimal order statistic
$X_{n:1} \equiv \min \{ X_1 , \ldots , X_n \}$.  By
\citet{MR0008655} 
or \citet[Theorem 1.1.2, page 5]{MR2234156}), it is well-known that
there exists a sequence $\{a_{n}\}$ such that 
\begin{eqnarray}
a_{n}^{-1} X_{n:1} \rightarrow_d Y 
\label{ConvMinOrderStat}
\end{eqnarray} 
where $Y$ has a
nondegenerate limiting distribution $G$ if and only if
\begin{eqnarray}
n F_0 (a_n x) \rightarrow x^{\gamma},  \qquad x > 0, 
\label{SequentialConvergenceRightTailCond}
\end{eqnarray}
for some $\gamma > 0$, 
and hence $a_n \rightarrow 0$.  
One possible choice of $a_n $ is $a_n = F_0^{-1} (1/n)$, but any sequence
$\{ a_n \}$ satisfying $nF_0 (a_n ) \rightarrow 1$ also works.
Since $F_{0}$ is concave the convergence in 
(\ref{SequentialConvergenceRightTailCond}) is uniform on any interval $[0,K]$.
Concavity of $F_0$ and existence of $f_0$ also implies
convergence of the derivative:
\begin{eqnarray}
&na_{n}f_{0}(a_{n}x)\to \gamma x^{\gamma-1}.
\label{ConseqGnedenkoHypothesisForDensity}
\end{eqnarray}
By \citet{MR0008655}, 
(\ref{ConvMinOrderStat}) is equivalent to
\begin{eqnarray}
\lim_{x\to 0+} \frac{F_{0}(cx)}{F_{0}(x)} = c^{\gamma},  \qquad c > 0.
\label{GnedenkoHypothesis}
\end{eqnarray}
Thus (\ref{ConvMinOrderStat}), (\ref{SequentialConvergenceRightTailCond}), and 
(\ref{GnedenkoHypothesis}) are equivalent.
In this case we have:
\begin{eqnarray}
G(x) = 1 - e^{-x^{\gamma}},\quad x\ge 0.
\label{GnedenkoLimitDistribution}
\end{eqnarray}
Since $F_{0}$ is concave, the power
$\gamma \in (0, 1]$.

As illustrations of our general result, we consider the  following
three hypotheses on $f_0$:
\begin{list}{}
        {\setlength{\topsep}{3pt}
        \setlength{\parskip}{0pt}
        \setlength{\partopsep}{0pt}
        \setlength{\parsep}{0pt}
        \setlength{\itemsep}{1pt}
        \setlength{\leftmargin}{30pt}}
\item[\textbf{G0:}]
The density $f_0$ is bounded at zero:  $f_0 (0) < \infty$.
\item[\textbf{G1:}]
For some $\beta \ge 0$ and $0 < C_1 < \infty$,
$$
(\log (1/x))^{-\beta} f_0(x) \rightarrow C_1 \qquad \mbox{as} \ \ x \searrow 0 .
$$
\item[\textbf{G2:}]
For some $0 \le \alpha < 1$ and $0 < C_2 < \infty$
$$
x^{\alpha} f_0(x) \rightarrow C_2 \qquad \mbox{as} \ \ x \searrow 0 .
$$
\end{list}
Note that in {\bf G2} the value $\alpha =1$ is not possible for a positive limit $C_2$ since $x f(x) \rightarrow 0$
as $x \rightarrow 0$
for any monotone density $f$; see e.g. \citet[Theorem 6.2, page 173]{MR836973}.
Below we assume that $F_{0}$ satisfies the condition
(\ref{GnedenkoHypothesis}).
Our cases $\mathbf{G0}$ and $\mathbf{G1}$ correspond to
$\gamma = 1$ and $\mathbf{G2}$ to $\gamma = 1-\alpha$.

One motivation for considering monotone densities which are unbounded at zero comes from the study of mixture
models.  An example of this type, as discussed by \citet{MR2065195}, is as follows.
Suppose $X_1, \ldots , X_n$ are i.i.d. with distribution function $F$ where,
\begin{align*}
&\hspace{0cm} \mbox{under} \ H_0 : \ F = \Phi,  \hspace{1cm} \mbox{the standard normal d.f.} \\
&\hspace{0cm}  \mbox{under} \ H_1 : \ F = (1-\epsilon) \Phi + \epsilon \Phi (\cdot - \mu ) ,  \ \ \epsilon \in (0,1),  \ \ \mu > 0.
\end{align*}
If we transform to $Y_i \equiv 1-\Phi (X_i) \sim G$, then, for $0 \le y \le 1$,
\begin{align*}
&\hspace{1cm}  \mbox{under} \ H_0 : G(y) = y, \hspace{0.75cm} \mbox{the Uniform}(0,1) \ \mbox{d.f.}, \\
&\hspace{1cm}  \mbox{under} \ H_1 :  G = G_{\epsilon, \mu} (y) = (1-\epsilon) y + \epsilon (1- \Phi ( \Phi^{-1}(1-y) - \mu )).
\end{align*}
It is easily seen that the density $g_{\epsilon, \mu}$ of $G_{\epsilon, \mu}$, given by
$$
g_{\epsilon, \mu} (y) = (1-\epsilon) + \epsilon \frac{\phi ( \Phi^{-1} (1-y) - \mu)}{\phi (\Phi^{-1} (1-y))},
$$
is monotone decreasing on $(0,1)$ and is unbounded at zero.  As we will show in
Section~\ref{sec:MixtureApplications},  $G_{\epsilon, \mu}$  satisfies
our key hypothesis (\ref{GnedenkoHypothesis}) below with $\gamma =1$.  Moreover, we will show
that the whole class of models of this type with $\Phi$ replaced by the generalized Gaussian
(or Subbotin) distribution, also
satisfy (\ref{GnedenkoHypothesis}), and hence the behavior of the Grenander
estimator at zero gives information about the behavior of the contaminating component of the mixture model (in the
transformed form) at zero.

Another motivation for studying these questions in the monotone density framework is to gain
insights for a study of the corresponding questions in the context of nonparametric estimation of
a monotone spectral density.  In that (related, but different) setting, singularities at the origin correspond
to the interesting phenomena of long-range dependence and long-memory processes; see e.g.
\citet{Cox:long:1984},  
\citet{MR1304490},     
\citet{MR1464601},     
\citet{MR1808873},     
and
\citet{MR1908951}.     
Although our results here do not apply directly to the problem of nonparametric estimation of
a monotone spectral density function, it seems plausible that similar results will hold in that setting;
note that when $f$ is a spectral density,  the assumptions {\bf G1} and {\bf G2} correspond to
long-memory processes (with the usual description being in terms of $\beta = 1-\alpha \in (0,1)$
or the Hurst coefficient $H = 1 - \beta /2 = 1 - (1-\alpha)/2 = (1+ \alpha)/2$).
See \citet{Anevski-Soulier:09} for recent work on nonparametric estimation of a monotone spectral density.

Let $\NN $ denote the standard Poisson process on $\RR^+$.
When (\ref{GnedenkoHypothesis}) and hence also (\ref{GnedenkoLimitDistribution}) hold,
it follows from
\citet[Theorem 2.1, page 522]{MR0415796} together with \citet[Theorem 2.15(c)(ii), pages 306-307]{MR1943877}, that
\begin{eqnarray}
n\mathbb{F}_{n}(a_{n}t)\Rightarrow \NN(t^{\gamma})\qquad \text{ in }D[0,\infty),
\label{GeneralWeakConvergenceECPGnedenkoHypothesis}
\end{eqnarray}
which should be compared to (\ref{SequentialConvergenceRightTailCond}).

Since we are studying the estimator $\widehat{f}_n$ near zero and because the value of $\widehat{f}_n$ at
zero is defined as the right limit $\lim_{x\searrow 0} \widehat{f}_n (x) \equiv \widehat{f}_n (0)$, it is sensible
to study instead the right-continuous modification of $\widehat{f}_n$, and this of course coincides with the
right derivative $\widehat{f}_n^R$ of the least concave majorant $\widehat{F}_n $ of the empirical distribution
function $\FF_n$.  
Therefore we change notation for
the rest of this paper and write   $\widehat{f}_n$ for $\widehat{f}_n^R$ throughout the following.
We write $\widehat{f}_n^{L}$
for the left-continuous Grenander estimator.

We now obtain the following theorem concerning the behavior of the Grenander estimator at zero.

\begin{thm}
\label{thm:GeneralTheoremGrenanderAtZero}
Suppose that (\ref{GnedenkoHypothesis}) holds.  Let $a_n$ satisfy $n F_0 (a_n) \sim 1$, 
let  $ \widehat{h}_{\gamma}$
denote the right derivative of the least concave majorant of $t \mapsto \NN(t^{\gamma})$, $t \ge 0$.
Then:\\
(i) \
$n a_n \widehat{f}_n (t a_n ) \Rightarrow \widehat{h}_{\gamma} (t) \qquad \mbox{in} \ \ D[0,\infty)$.\\
(ii) \ For all $c \ge 0$
\begin{eqnarray*}
\sup_{0<x\le c a_n} \left|\frac{\widehat{f}_{n}(x)}{f_{0}(x)} - 1\right|
\to_{d} \sup_{0<t\le c}\left|\frac{t^{1-\gamma}\widehat{h}_{\gamma} (t)}{\gamma} - 1\right |.
\end{eqnarray*}
\end{thm}
\medskip

The behavior of $\widehat{f}_n$ near zero under the different hypotheses {\bf G0}, {\bf G1}, and {\bf G2}
now follows as corollaries to Theorem 1.1.
Let $Y_{\gamma} \equiv \widehat{h}_{\gamma} (0)$.  We then have
\begin{eqnarray}
Y_{\gamma} = \sup_{t>0} (\NN (t^{\gamma})/t) = \sup_{s>0} ( \NN (s) / s^{1/\gamma} ).
\label{DefinitionOfLimitingVariableAtZero}
\end{eqnarray}
Here we note that $Y_1 =_d 1/U$ where $U \sim \mbox{Uniform}(0,1)$ has
distribution function $H_1 (x) = 1 - 1/x$ for  $x \ge 1$.
The distribution of $Y_{\gamma}$ for $\gamma \in (0,1]$
 is given in Proposition~\ref{DistributionOfYSubAlpha} below.
The first part of the following corollary was established by \citet{MR1243398}.

\begin{cor}
\label{cor:BoundedCase}
Suppose that {\bf G0} holds.  Then $\gamma =1$, $a_n^{-1} = n f_0(0+)$ 
satisfies $n F_0 (a_n) \rightarrow 1$, and it follows that: \\
(i) \
$$
\widehat{f}_n (0) \rightarrow_d f_0 (0) \widehat{h}_1 (0)
     = f_0 (0) Y_1 .
$$
(ii) \  The processes $\{ t \mapsto \widehat{f}_n (t n^{-1} ) : \ n \ge 1 \}$ satisfy
$$
\widehat{f}_n (t n^{-1} ) \Rightarrow f_0 (0) \widehat{h}_1 (f_0 (0) t)  \qquad \mbox{in} \ \ D[0, \infty).
$$
(iii) \ For $c_n = c/n$ with $c> 0$,
\begin{eqnarray*}
\sup_{0< x \leq c_n} \left|\frac{\widehat{f}_n(x)}{f_0(x)}-1\right| \rightarrow_d Y_1 -1
\end{eqnarray*}
which has distribution function $H_1(x+1)=1-1/(x+1)$ for $x\geq 0.$
\end{cor}

\begin{cor}
\label{cor:LogGrowthCase}
Suppose that {\bf G1} holds.
Then $F_0 (x) \sim C_1 x ( \log (1/x) )^{\beta}$, so  $\gamma =1$, and 
$a_n^{-1} = C_1 n (\log n)^{\beta}$ satisfies $nF_0 (a_n) \rightarrow 1$.  It follows that:\\
(i) \
\begin{eqnarray*}
\frac{\widehat{f}_n (0)}{(\log n)^{\beta} } \rightarrow_d C_1 Y_1  .
\end{eqnarray*}
(ii) \ The processes $\{ t \mapsto (\log n)^{-\beta} \widehat{f}_n (t /( n (\log n)^{\beta})  ) : \ n \ge 1 \}$
satisfy
\begin{eqnarray*}
\frac{1}{(\log n)^{\beta}} \widehat{f}_n \left ( \frac{t}{n (\log n)^{\beta} } \right )
     \Rightarrow C_1 \widehat{h}_1 (C_1 t) \qquad \mbox{in} \ \ D[0,\infty)
\end{eqnarray*}
(iii) \ For $c_n = c / (n (\log n)^{\beta} )$ with $c>0$,
\begin{eqnarray*}
\sup_{0< x \leq c_n } \left|\frac{\widehat{f}_n(x)}{f_0(x)}-1\right|\rightarrow_d Y_1-1.
\end{eqnarray*}
\end{cor}
\medskip

\begin{cor}
\label{cor:PolyGrowthCase}
Suppose that {\bf G2} holds and set $\widetilde{C}_2 = (C_2/(1-\alpha))^{1/(1-\alpha)}$.
Then $F_0 (x) \sim C_2 x^{1-\alpha}/(1-\alpha)$, so  $\gamma =1-\alpha$,
$a_n^{-1} = \widetilde{C}_2 n^{1/(1-\alpha)}$ 
satisfies $nF_0 (a_n) \rightarrow 1$, and it follows that:\\
(i) \
\begin{eqnarray}
\frac{\widehat{f}_n (0)}{n^{\alpha/(1-\alpha)}}
 \rightarrow_d  \widetilde{C}_2 Y_{1-\alpha}  .
\label{LimitDistrib}
\end{eqnarray}
(ii) \ The processes
$\{ t \mapsto n^{-\alpha/(1-\alpha)} \widehat{f}_n ( t n^{-1/(1-\alpha)} ): \ n \ge 1 \}$
satisfy
\begin{eqnarray*}
\frac{\widehat{f}_n ( t n^{-1/(1-\alpha)} )}{n^{\alpha/(1-\alpha)}} \Rightarrow
    \widetilde{C}_2 \widehat{h}_{1-\alpha} ( \widetilde{C}_2   t  )
\qquad \mbox{in} \ \ D[0, \infty ) .
\end{eqnarray*}
(iii) \ For $c_n = c/ n^{1/(1-\alpha)}$ with $c>0$,
\begin{eqnarray*}
\sup_{0< x \leq c_n} \left|\frac{\widehat f_n(x)}{f_0 (x)}-1\right|
  \rightarrow_d \sup_{0< t \leq c \widetilde{C}_2} \left|\frac{t^\alpha \widehat{h}_{1-\alpha} (t)}{1-\alpha}-1\right| .
\end{eqnarray*}
\end{cor}

Taking $\beta = 0 $ in (i) of Corollary~\ref{cor:LogGrowthCase} yields
the limit theorem (\ref{WoodroofeSunThmBoundedCase}) of
\citet{MR1243398}  
as a corollary; in this case $C_1= f_0(0)$.
Similarly,  taking $\alpha = 0 $ in (ii) of  Corollary~\ref{cor:PolyGrowthCase} yields
the limit theorem (\ref{WoodroofeSunThmBoundedCase}) of
\citet{MR1243398}  
as a corollary; in this case $C_2= f_0(0)$.
Note that Theorem~\ref{thm:GeneralTheoremGrenanderAtZero}
yields further corollaries when assumptions
{\bf G1} and {\bf G2} are modified by other
slowly varying functions.

Recall the definition (\ref{DefinitionOfLimitingVariableAtZero}) of $Y_{\gamma}$.
The following proposition gives the distribution of $Y_{\gamma}$ for $\gamma \in (0,1]$.

 \begin{prop}
 \label{DistributionOfYSubAlpha}
        For fixed $0 < \gamma \le  1$ and $x>0$,
        \[
            \Pr\left( \sup_{s>0}\left\{
            \frac{\NN(s)}{s^{1/\gamma}}\right\}
            \leq x \right) = \left \{
            \begin{array}{l l}
                  1- 1/x \, , & \mbox{if} \ \gamma = 1, \ x\ge 1,\\
                  1 - \sum_{k=1}^{\infty} a_k(x,\gamma)\,, &  \mbox{if} \ \gamma < 1 , \ x > 0 ,
            \end{array} \right .
        \]
        where the sequence $\{ a_k(x,\gamma) \}_{k\geq 1}$ is
        constructed recursively as follows:
        \begin{align*}
            a_1(x,\gamma) &= p\left( \left( \frac{1}{x} \right)^{\gamma} ; 1\right)\,,\\
                \intertext{and, for $j\geq 1$,}
            a_k(x,\gamma) &= p\left( \left( \frac{k}{x} \right)^{\gamma}
                                           ; k \right) - \sum_{i=1}^{k-1}\left\{ a_i(x,\gamma)\cdot
            p\left( \left( \frac{k}{x} \right)^{\gamma}
                              - \left(\frac{i}{x} \right)^{\gamma}
                                ; k-i\right)\right\}\,,
        \end{align*}
        where $p(m; k) \equiv e^{-m}m^k/k!$.
    \end{prop}

\begin{figure}[htb!]
\centering
\includegraphics[width=0.8\textwidth, height=0.5\textwidth]{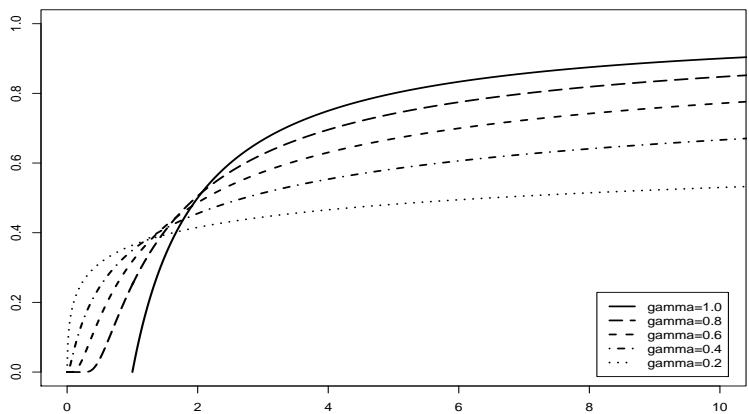}
\caption{The distribution functions of $Y_{\gamma}$, $\gamma \in \{0.2, 0.4, 0.6, 0.8, 1.0 \}$.}
\label{fig:DisributionsOfYSubGamma}
\end{figure}

\begin{rem}
\label{YSubAlphaHasHeavyTails}
The random variables $Y_{\gamma}$ are increasingly heavy-tailed as $\gamma $ decreases; cf.
Figure~\ref{fig:DisributionsOfYSubGamma}.
Let $T_1, T_2, \ldots $ be the event times of the Poisson process $\NN$; i.e. 
$\NN(t) = \sum_{j=1}^\infty 1_{[T_j \le t]}$.  Then
note that
\begin{eqnarray*}
Y_{\gamma} \stackrel{d}{=} \sup_{j\ge 1} \frac{j}{T_j^{1/\gamma}} \ge \frac{1}{T_1^{1/\gamma}}
\end{eqnarray*}
where $T_1 \sim \ $Exponential$(1)$.  On the other hand
\begin{eqnarray*}
Y_{\gamma}
 =  \left ( \sup_{t>0} \frac{\NN (t)^{\gamma}}{t} \right )^{1/\gamma}
 \le  \left ( \sup_{t>0} \frac{\NN(t)}{t} \right )^{1/\gamma} \stackrel{d}{=} \frac{1}{U^{1/\gamma}}
\end{eqnarray*}
where $U \sim \ $Uniform$(0,1)$.    Thus it is easily seen that $E( Y_{\gamma}^r ) < \infty$ if and only if $r < \gamma$,
and that the distribution function $F_{\gamma}$ of $Y_{\gamma}$ is bounded above and below by the
distribution functions $G_{\gamma}^L$ and $G_{\gamma}^U$ of $ 1/T_1^{1/\gamma}$ and $1/U^{1/\gamma}$, respectively.
\end{rem}

The proofs of the above results appear in  Appendix A.  
They rely heavily on a set equality known as the ``switching relation".  
We study this relation using convex analysis in Section \ref{sec:switching}.  
Section \ref{sec:Numerical} gives some numerical results which accompany 
the results presented here, and Section \ref{sec:MixtureApplications} studies 
applications to the estimation of mixture models.

\section{Switching relations}
\label{sec:switching}

In this section we consider several general variants of the so-called
switching relation first given in
\citet{MR822052},   
and used repeatedly by other authors, including
\citet{MR2211085, MR2283391},
and
\citet{vanderVaartWellner:96}.
Other versions of the switching relation were also studied by
\citet[Lemma 4.1]{vanderVaartvanderLaan:06}.
In particular, we provide a novel proof of the result using convex analysis.
This approach also allows us to re-state the relation without restricting the domain to compact intervals.
Throughout this section we make use of definitions from convex analysis 
(cf. \citet{MR0274683, MR1491362, MR2061575}) which are given in Appendix B.

Suppose that $\Phi$ is a function, $\Phi : D \rightarrow \RR$, defined on the (possibly infinite) closed
interval $D \subset \RR$.
The least concave majorant $\widehat{\Phi}$ of $\Phi$ is the pointwise infimum of all
closed concave functions
$g: D \rightarrow \RR$ with $g \ge \Phi$.
Since $\widehat{\Phi}$ is concave, it is continuous on $D^o$, the interior of $D$.
Furthermore, $\widehat{\Phi}$ has left and right derivatives on $D^o$, and is differentiable with the
exception of at most countably many points.
Let $\widehat{\phi}_L$ and $\widehat{\phi}_R$ denote the left and right derivatives, respectively,
of $\widehat{\Phi}$.

If $\Phi$ is upper semicontinuous, then so is the function $\Phi_y (x)= \Phi (x) - y x$ for each $y \in \RR$.
If $D$ is compact, then $\Phi_y$ attains a maximum on $D$, and the set of points achieving the maximum
is closed.  Compactness of $D$ was assumed by
\citet[see their Lemma 4.1, page 24]{vanderVaartvanderLaan:06}.
One of our goals here is to relax this assumption.

Assuming they are defined, we consider the argmax functions
\begin{eqnarray*}
 \kappa_L (y) \equiv \mbox{argmax}^{L} \Phi_y
    & \equiv & \mbox{argmax}_x^L \{ \Phi (x) - yx \} \\
    & = & \inf \{ x \in D: \ \Phi_y (x) = \sup_{z \in D} \Phi_y (z) \},\\
\kappa_R (y) \equiv \mbox{argmax}^{R} \Phi_y
     & \equiv & \mbox{argmax}_x^R \{ \Phi (x) - yx \} \\
     & = & \sup\{ x \in D : \ \Phi_y (x) = \sup_{z \in D} \Phi_y (z) \}.
\end{eqnarray*}

\begin{thm}
\label{GeneralSwitchingRelations}
Suppose that $\Phi$ is a proper upper-semicontinuous real-valued function defined on a closed subset
$D \subset \RR$.
Then $\widehat{\Phi}$ is proper if and only if $\Phi \le l$ for some linear function $l$ on $D$.
Furthermore, if $\mbox{conv} ( \mbox{hypo} (\Phi))$ is closed, then the functions $\kappa_L $ and
$\kappa_R$ are well defined and the following two switching relations hold:
 for $x \in D$ and $y \in \RR$,
\begin{description}
\item[S1]
$\widehat{\phi}_L (x) < y$ if and only if $\kappa_R (y) < x$.
\item[S2]
$\widehat{\phi}_R (x) \le y$ if and only if $\kappa_L (y) \le  x$.
\end{description}
\end{thm}

When $\Phi$ is the empirical distribution function $\FF_n$  as in
Section~\ref{sec:IntroMainResults}, then
$\widehat{\Phi} = \widehat{F}_n$ is the least concave majorant of $\FF_n$,
and $\widehat{\phi}_L = \widehat{f}_n^L$ the Grenander estimator as defined in
Section~\ref{sec:IntroMainResults},
while  $\widehat{\phi}_R = \widehat{f}_n = \widehat{f}_n^R$ is the right continuous version of the estimator.
In this situation the argmax functions $\kappa_R,\kappa_L$ correspond to
\begin{eqnarray*}
\widehat{s}^R_n (y)
& = & \sup \{ x \ge 0 : \FF_n(x) - yx  = \sup_{z \ge 0} (\FF_n (z) - yz ) \}, \\
\widehat{s}^L_n (y)
& = & \inf \{ x \ge 0 : \FF_n(x) - yx  = \sup_{z \ge 0} (\FF_n (z) - yz ) \}.
\end{eqnarray*}
The switching relation given by \citet{MR822052} says that with probability one
\begin{eqnarray}
\{ \widehat{f}_n^L (x) \le y\} = \{ \widehat{s}_n^R (y) \le x \} .
\label{SwitchingRelationClaimed}
\end{eqnarray}
\citet[page 296]{vanderVaartWellner:96}, say that (\ref{SwitchingRelationClaimed})
holds for every $x$ and $y$; see also
\citet[page 2229]{MR2211085}, and \citet[page 744]{MR2283391}.
The advantage of (\ref{SwitchingRelationClaimed}) is immediate:  the MLE is related to a
continuous map of a process whose behavior is well-understood.

The following corollary gives the conclusion of Theorem~\ref{GeneralSwitchingRelations}
when $\Phi$ is the empirical distribution function $\FF_n$.

\begin{cor}
\label{TrueSwitchingRelationsForGrenander}
Let $\widehat{F}_n$ be the least concave majorant of the empirical distribution function
$\FF_n$, and let $\widehat{f}_n^L$ and $\widehat{f}_n^R$ denote its left and right derivatives
respectively.  Then:
\begin{eqnarray}
&& \{ \widehat{f}_n^L (x) < y \} = \{ \widehat{s}_n^R (y) < x \},
       \label{OpenSwitchingRelationGrenander}    \\
&& \{ \widehat{f}_n^R (x) \le y \} = \{ \widehat{s}_n^L (y) \le x \} .
\label{ClosedSwitchingRelationGrenander}
\end{eqnarray}
\end{cor}

The following example  shows, however,
that the set identity (\ref{SwitchingRelationClaimed}) can fail.

\begin{expl}
Suppose that we observe $(X_1 , X_2 , X_3) = (1,2,4)$.  Then the MLE $\widehat{f}_n^L$
is given by
$$
\widehat{f}_n^L (x) = \left \{ \begin{array}{l l} 1/3, & 0 < x \le 2, \\
                                                                  1/6, & 2 < x \le 4, \\
                                                                  0, & 4 < x < \infty .
                                       \end{array} \right .
$$
The process $\widehat{s}_n^R $ is given by
$$
\widehat{s}_n^R (y) = \left \{ \begin{array}{l l} 4, & 0 < y \le 1/6, \\
                                                                   2, & 1/6 < y \le 1/3, \\
                                                                  0, & 1/3 < y < \infty .
                                       \end{array} \right .
$$
Note that (\ref{SwitchingRelationClaimed}) fails if $x=4$ and $0 < y < 1/6$, since in this case
$\widehat{f}_n^L (x) =\widehat{f}_n^L (4) = 1/6$ and the event $\{ \widehat{f}_n^L (x) \le y \}$ fails to hold
while $\widehat{s}_n^R (y) = 4$ and the event $\{ \widehat{s}_n^R (y) \le x \}$ holds.
However, (\ref{OpenSwitchingRelationGrenander}) does hold:  with $x=4$ and $0 < y < 1/6$,
both of the events  $\{ \widehat{f}_n^L (x) < y \}$ and $\{ \widehat{s}_n^R (y) < x \}$ fail to hold.
Some checking shows that (\ref{OpenSwitchingRelationGrenander}) as well as
(\ref{ClosedSwitchingRelationGrenander})
hold for all other values
of $x$ and $y$.
\end{expl}

Our proof of Theorem~\ref{GeneralSwitchingRelations} will be based on the following
proposition which is a consequence of general facts concerning convex functions as given
in \citet{MR0274683} and \citet{MR1491362}. 

\begin{prop}
\label{ConvexDualityProp}
Let $h$ be a closed proper convex function on $\RR$, and let $f$ be its conjugate,
$$
f(y) = \sup_{x \in \RR} \{ y x - h(x) \} .
$$
Let $h_{-}'$ and $h_{+}'$ be the left and right derivatives of $h$, and define functions $s_{-}$ and $s_{+}$ by
\begin{eqnarray}
&& s_{-} (y) = \inf \{ x \in \RR : \ yx - h(x) = f(y) \},
\label{LeftVersionArgmax}\\
&& s_{+} (y) = \sup\{ x \in \RR : \ yx - h(x) = f(y) \} .
\label{RightVersionArgmax}
\end{eqnarray}
Then the following set identities hold:
\begin{eqnarray}
&& \{ (x,y) : \ h'_{-} (x) \le y \} = \{ (x,y) : \ s_+ (y) \ge x \} , \label{SetIdentityOne} \\
&& \{ (x,y) : \ h'_{+} (x) < y \} = \{ (x,y) : \ s_-(y) >  x \} , \label{SetIdentityTwo}
\end{eqnarray}
\end{prop}
\medskip

\par\noindent
{\bf Proof.}
All the references in this proof are to \citet{MR0274683}.
By Theorem 24.3 (page 232) the set
$\Gamma = \{ (x,y) \in \RR^2 : \ y \in \partial h(x) \}$ (i.e. the graph of $\partial h$),
is a maximal complete non-decreasing curve. By Theorem 23.5, page 218,  the closed proper convex function $h$ and its conjugate
$f$ satisfy
$$
h(x) + f(y) \ge xy
$$
and equality holds if and only if $y \in \partial h (x)$, or equivalently if
$x \in \partial f (y)$ where $\partial h$ and $\partial f$ denote the subdifferentials of
$h$ and $f$ respectively (see page 215). Thus we also have:
$$
\Gamma = \{ (x,y) \in \RR^2 : \ x \in \partial f(x) \},
$$
and, by the definitions of $s_{-}$ and $s_{+}$,
$$
\Gamma = \{ (x,y) : \ s_{-} (y) \le x \le s_{+} (y) \}.
$$ 
By Theorem 24.1 (page 227) the curve $\Gamma$
is defined by the left and right derivatives of $h$:
\begin{eqnarray}
&& \Gamma = \{ (x,y) : \ h_{-}' (x) \le y \le h_{+}' (x) \} .
\label{CurveRepresentation1}
\end{eqnarray}
Using the dual representation we obtain:
\begin{eqnarray}
&& \Gamma = \{ (x,y) : \ f_{-}' (y) \le x \le f_{+}' (y) \},
\label{CurveRepresentation2}
\end{eqnarray}
therefore $s_{-} \equiv f_{-}'$ and $s_{+} \equiv f_{+}' $.
Moreover,  the functions $h_{-}'$ and $f_{-}'$ are left-continuous,
the functions $h_{+}'$ and $f_{+}'$ are right continuous, and all
of these functions are nondecreasing.

From \eqref{CurveRepresentation1} and \eqref{CurveRepresentation2} it follows that:
\begin{eqnarray*}
 \{ h_{-}' (x) \le y \} = \{ f_{+}' (y) \ge x \},
\end{eqnarray*}
which implies (\ref{SetIdentityOne}).
Since the functions $h$ and $f$ are conjugate to each other, the relations between
them are symmetric.  Thus we have
$$
\{ f_{-}' (y) \le x \} = \{ h_{+}' (x) \ge y \},
$$
or equivalently
$$
\{ f_{-}' (y) > x \} = \{ h_{+}' (x) < y \},
$$
which implies (\ref{SetIdentityTwo}).
\hfill $\Box$
\bigskip

Before proving Theorem~\ref{GeneralSwitchingRelations} we need the
following two lemmas.

\begin{lem}
\label{lem:SuperLevelSetsDefined}
Let $S = \argmax_{D} \Phi$ and
$\widehat{S} = \argmax_{D} \widehat\Phi$
 be the maximal superlevel sets of $\Phi$ and $\widehat{\Phi}$.
 Then the set $\widehat{S}$ is defined if and only if the set $S$ is defined and in this case
 $\conv(S) \subseteq \widehat{S}$.
\end{lem}

\begin{lem}
\label{lem:CondForEqualitySuperLevelSets}
If $\conv(\hypo(\Phi))$ is a closed convex set then $\conv(S) = \widehat{S}$.
\end{lem}

\par\noindent
{\bf Proof of Lemma~\ref{lem:SuperLevelSetsDefined}:}
Since $\cl(\Phi)\le \widehat\Phi$ the set $S$ is defined if $\widehat S$ is defined.
On the other hand, if $S$ is defined then $\Phi$ is bounded from above on $D$.
Since:
\[\sup_{D} \Phi = \sup_{D}\widehat\Phi,\]
the function $\widehat\Phi$ is also bounded from above on $D$,  i.e. the set $\widehat S$ is defined.

By~\eqref{sw2-hypo-embed1} we have $S\subseteq\widehat{S}$.
Since $\Phi$ and $\widehat\Phi$ are
upper semicontinuous the sets $S$ and $\widehat S$  are closed. Since $\widehat S$
is convex we have $\conv(S) \subseteq \widehat S$.
\hfill$\Box$
\smallskip

\par\noindent
{\bf Proof of Lemma~\ref{lem:CondForEqualitySuperLevelSets}:}
Indeed, we have $\conv(\hypo(\Phi))\equiv \conv(\cl(\hypo(\Phi)))$, and
\[
\conv(\hypo(\Phi))\subseteq \hypo(\widehat{\Phi}) .
\]
Therefore $\conv(\hypo(\Phi))$ is a hypograph of some closed concave function $H$ such that:
\[
\Phi \le H \le \widehat{\Phi}.
\]
Thus $H=\widehat{\Phi}$.
The set $\widehat{S}$ is a face of $\hypo(\widehat{\Phi})$ and the set $\conv(S)$ is a face of $\conv(\hypo(\Phi))$.
The statement now follows from \citet[Theorem 18.3, page 165]{MR0274683}.
\hfill $\Box$
\medskip

\par\noindent
{\bf Proof of Theorem~\ref{GeneralSwitchingRelations}.}
To prove the first statement, first suppose $\widehat{\Phi}$ is proper. We have:
\begin{align}
&\hypo(\Phi) \subseteq \hypo(\cl (\Phi)) \equiv \cl(\hypo(\Phi)) \subseteq \cl(\conv(\hypo(\Phi)))\equiv \hypo(\widehat\Phi)
\label{sw2-hypo-embed1}
\end{align}
and therefore $\hypo(\Phi)$ is bounded by any support plane of $\hypo(\widehat\Phi)$.
This implies that there exists a linear function $l$ such that $\Phi \le l$.

Now suppose that there exists a linear
function $l$ such that $\Phi\le l$ on $D$. Then $\cl(\Phi)\le l$ and from \eqref{sw2-hypo-embed1} we have:
\begin{align*}
&\hypo(\Phi)\subseteq\hypo(l), \\
&\conv(\hypo(\Phi))\subseteq\hypo(l), \\
&\hypo(\widehat{\Phi})  \equiv \cl(\conv(\hypo(\Phi)))\subseteq\hypo(l) .
\end{align*}
Thus $\widehat{\Phi} < +\infty$ on $D$.
Since $\hypo(\Phi) \subseteq \hypo(\widehat{\Phi})$ there exists a finite point in $\hypo(\widehat{\Phi})$.

To show that the two switching relations hold, first consider the
convex function $h = - \widehat{\Phi}$.  Then
\begin{eqnarray*}
&&\widehat{\phi}_{L} (x) = - h'_{-} (x) , \\
&&\widehat{\phi}_{R} (x) = - h'_{+} (x) , \\
&& \kappa_{L} (y) = s_{-} (-y), \\
&& \kappa_{R} (y) = s_{+} (-y) ,
\end{eqnarray*}
and by the properness of $\widehat{\Phi}$
proved above and
Proposition~\ref{ConvexDualityProp}, it suffices to show that
\begin{eqnarray*}
&& \mbox{argmax}_x^L ( \Phi (x) - yx) = \mbox{argmax}_x^L ( \widehat{\Phi} (x) - yx ) , \\
&& \mbox{argmax}_x^R ( \Phi (x) - yx) = \mbox{argmax}_x^R ( \widehat{\Phi} (x) - yx ) .
\end{eqnarray*}
To accomplish this, it suffices, without loss of generality, to prove the equalities in the last display
when $y=0$, and this in turn will follow if we relate the maximal superlevel sets of $\Phi$ and $\widehat{\Phi}$.
This follows from Lemmas~\ref{lem:SuperLevelSetsDefined}
and~\ref{lem:CondForEqualitySuperLevelSets}.
\hfill $\Box$

\begin{rem}
\label{CounterExample1}
Note that $\conv(S) \neq \widehat S$ in general.  To see this,
 consider the function $\Phi$ defined on $\RR$ as follows:
\[
\Phi(x) = \begin{cases}0 & x\neq 0\\ 1 & x=0.\end{cases}
\]
We have that $ \Phi$ is upper-semicontinuous,
$S = \{0\}$ and $\widehat{\Phi} \equiv 1$, so $\widehat{S} = \RR$.
\end{rem}

\begin{rem}
\label{SufficientCondForClosure}
Note that  if $\conv(\hypo(\Phi))$ is a polyhedral set, then it is
closed (see e.g. \citet[Corollary 19.1.2]{MR0274683}). This is the case in our applications.
\end{rem}
\medskip

\begin{figure}[htb!]
\centering
\includegraphics[width=0.49\textwidth, height=0.37\textwidth]{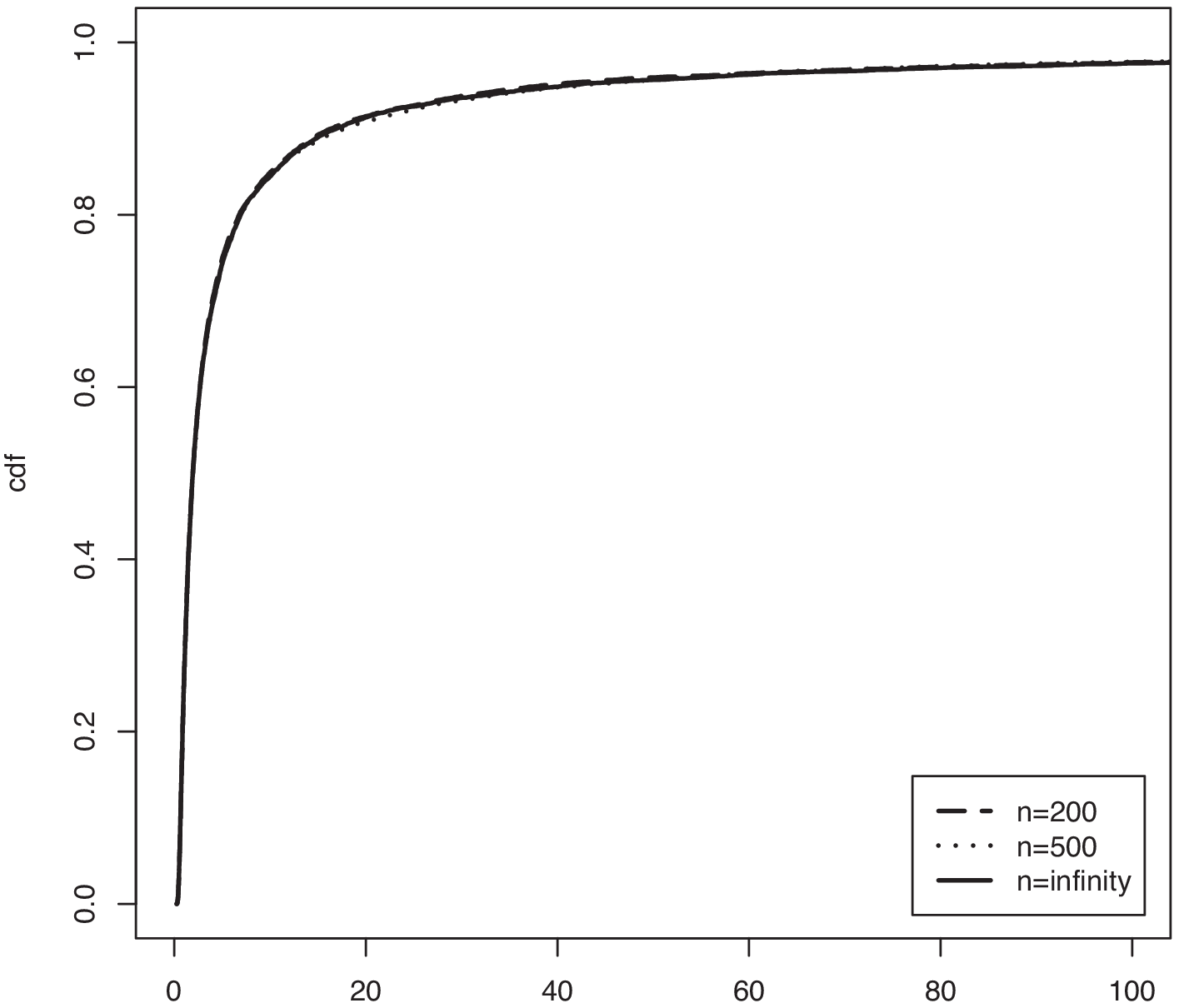}
\includegraphics[width=0.49\textwidth, height=0.37\textwidth]{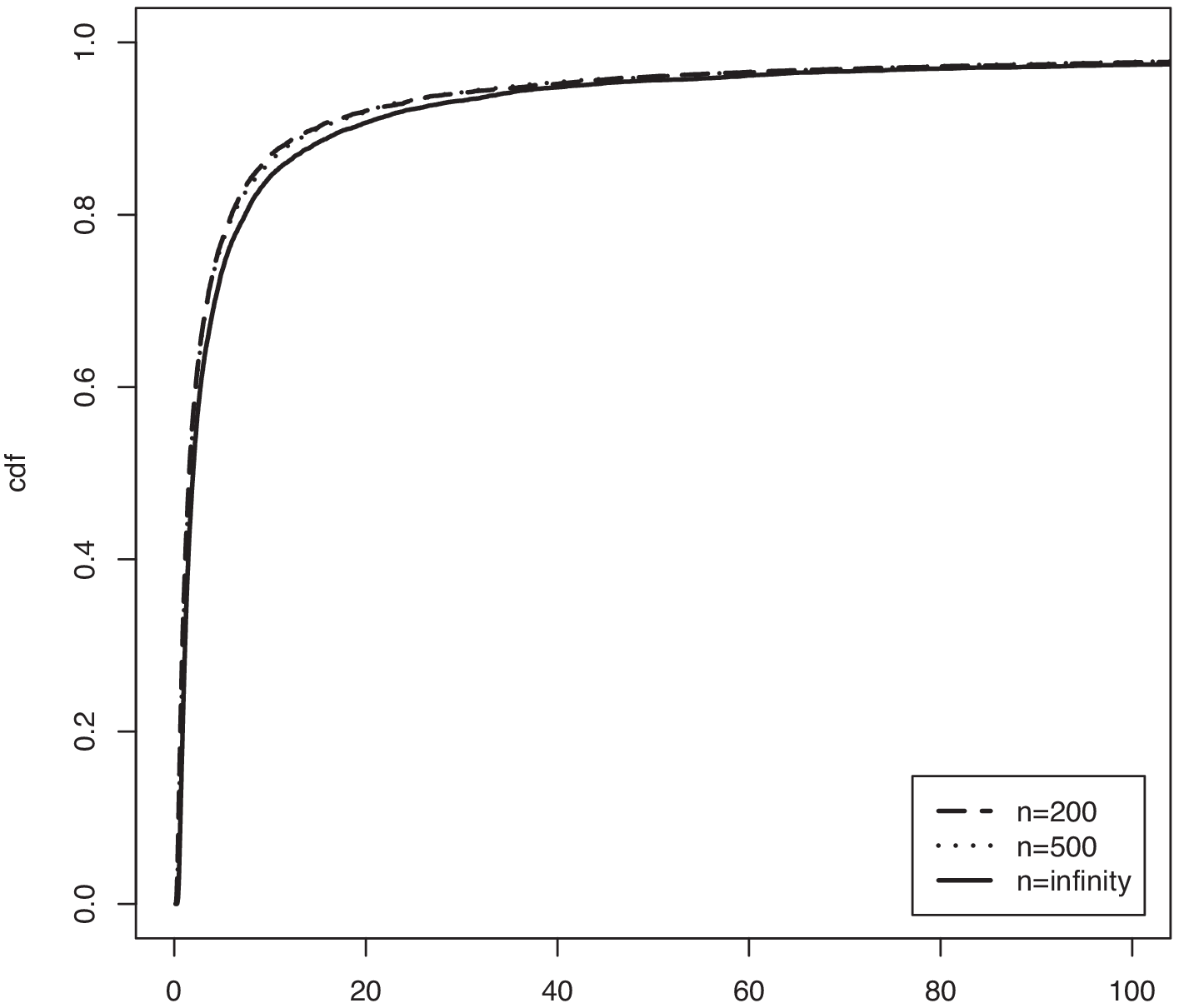}
\includegraphics[width=0.49\textwidth, height=0.37\textwidth]{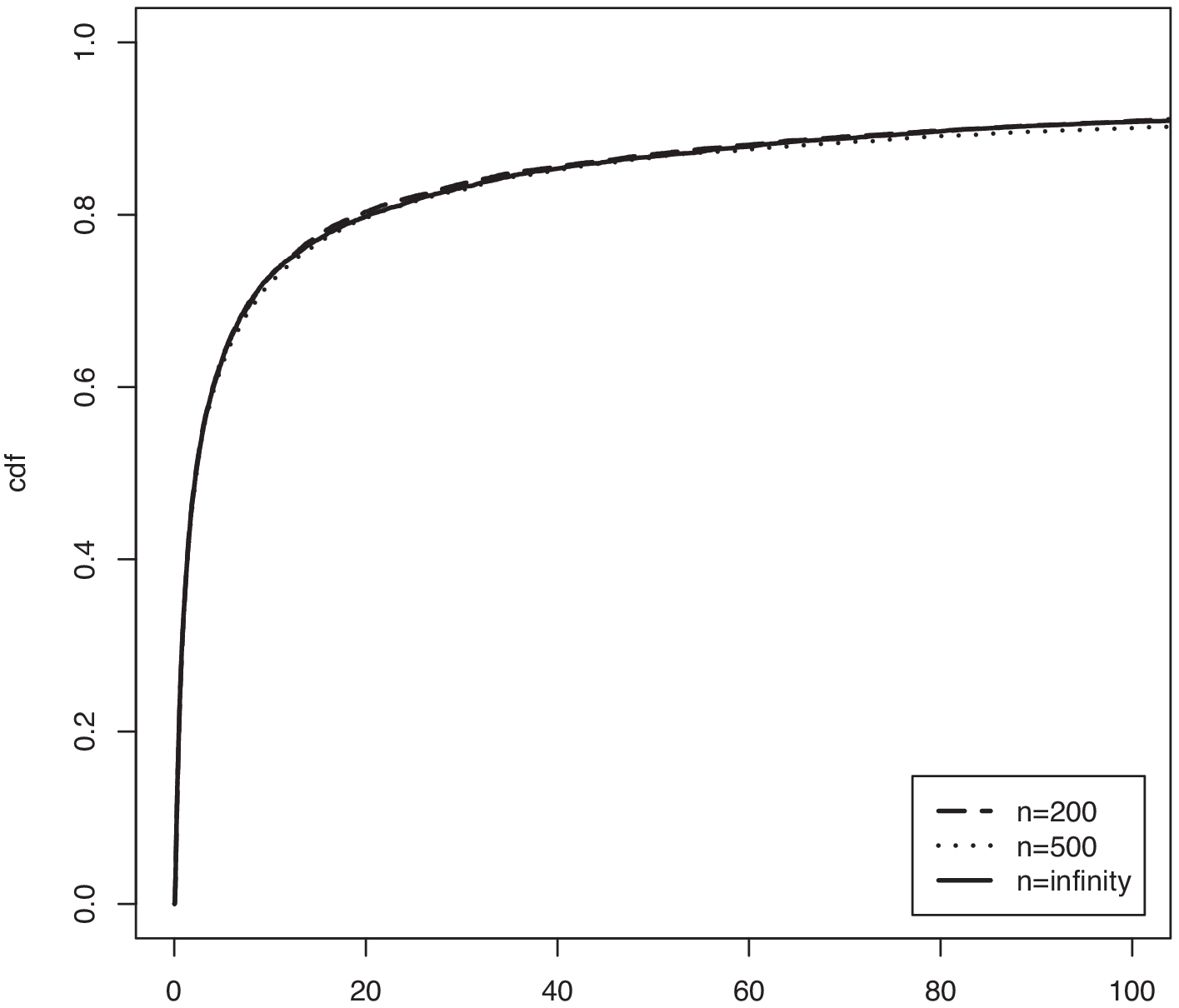}
\includegraphics[width=0.49\textwidth, height=0.37\textwidth]{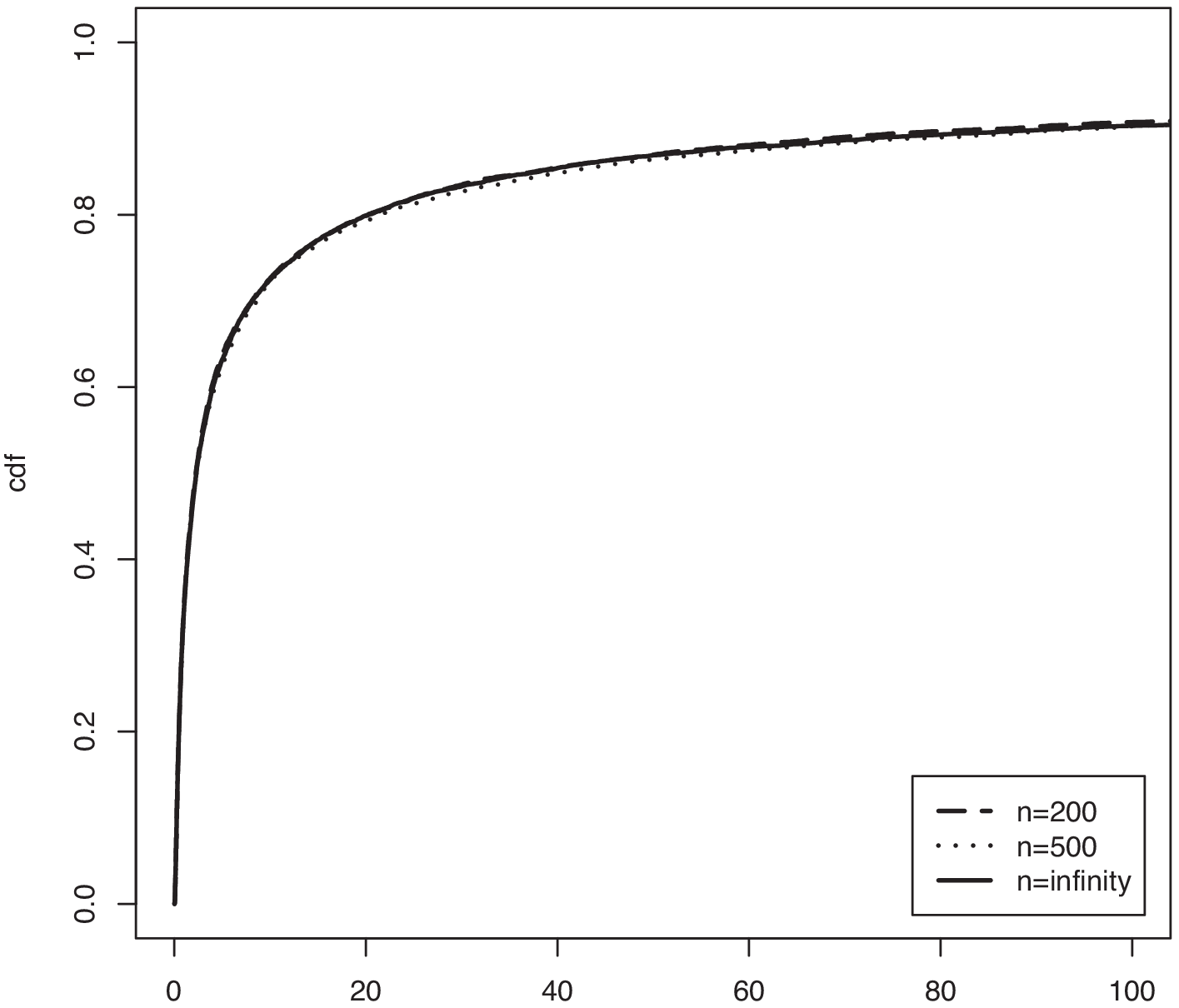}
\includegraphics[width=0.49\textwidth, height=0.37\textwidth]{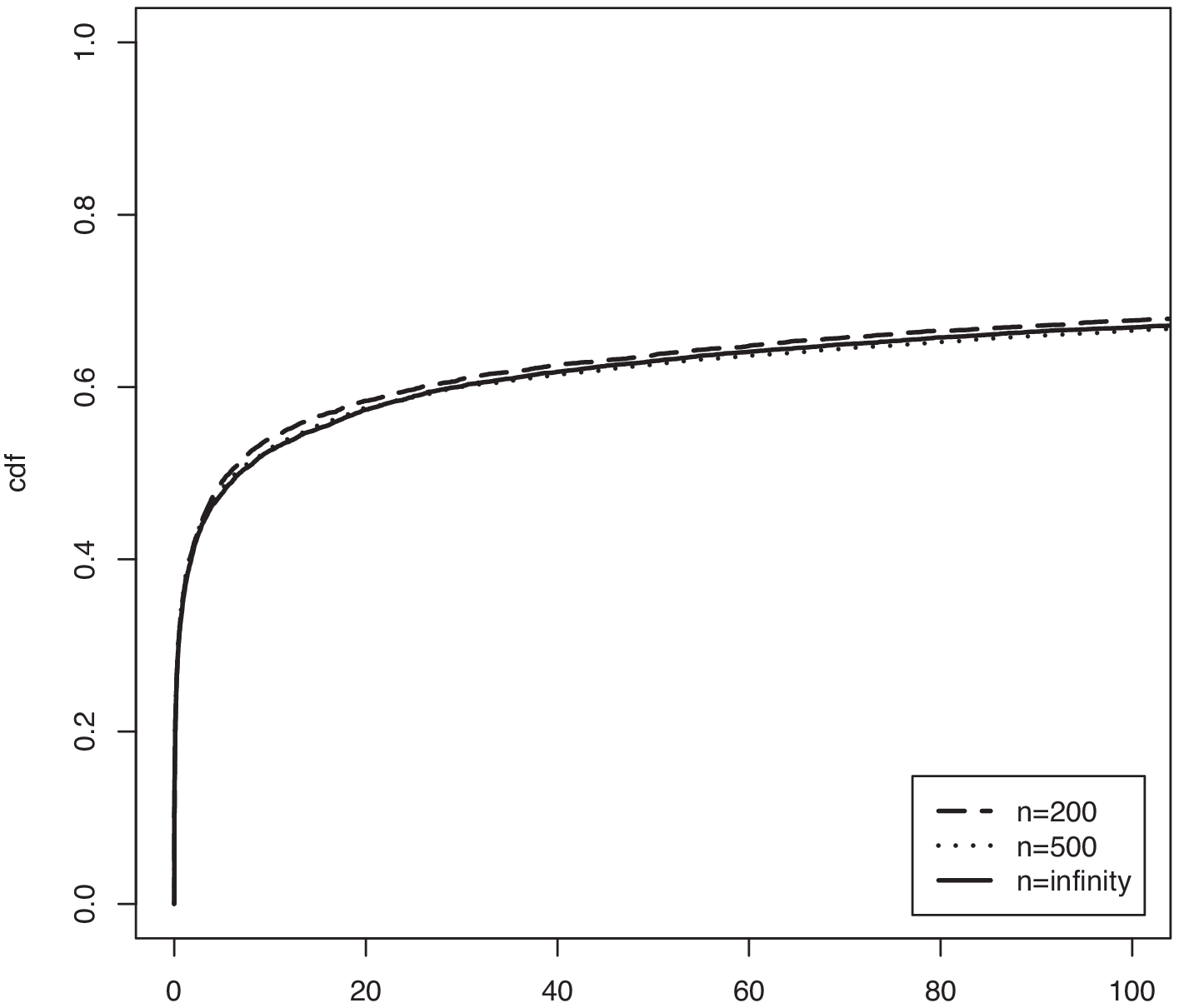}
\includegraphics[width=0.49\textwidth, height=0.37\textwidth]{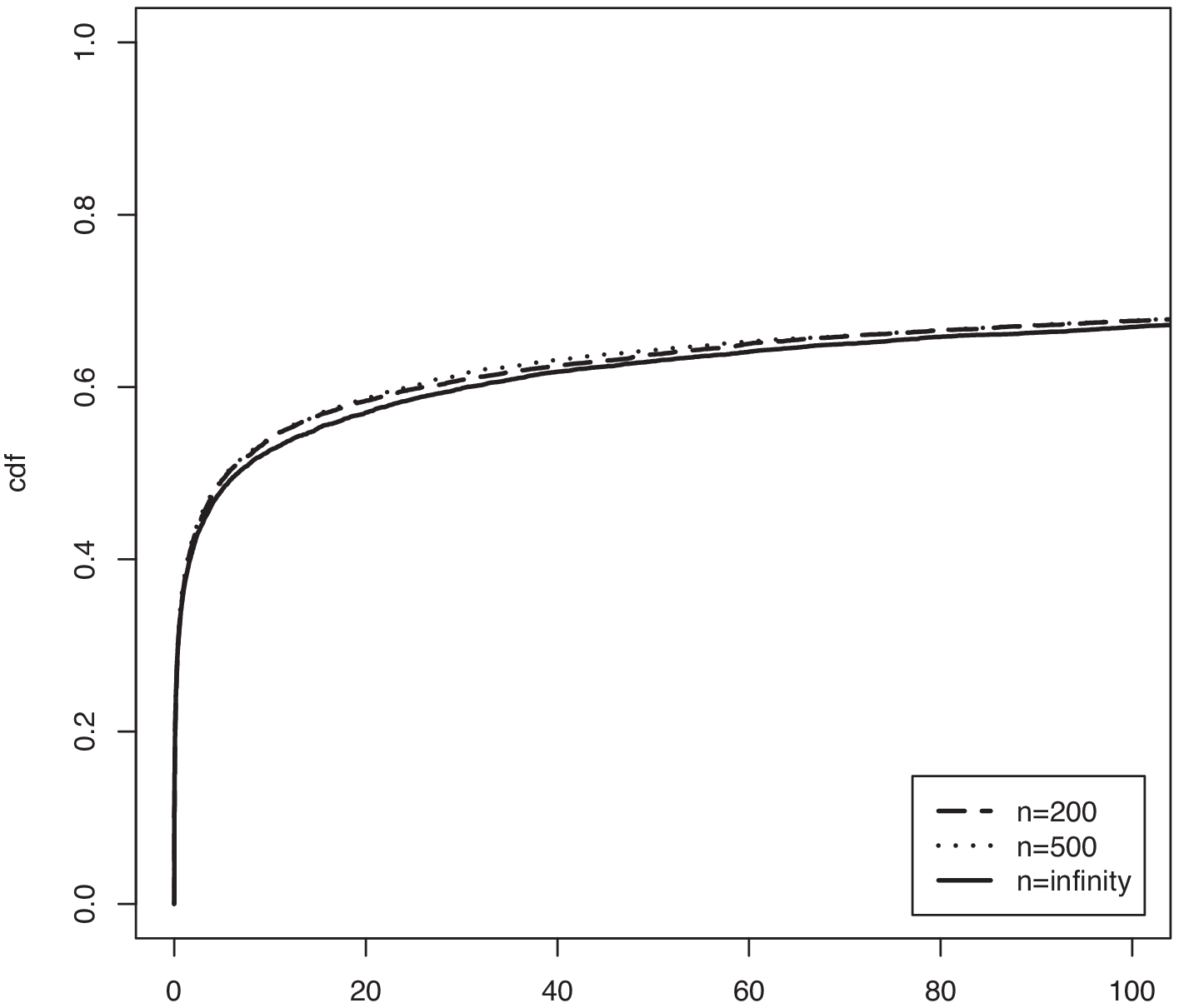} 
\caption{Empirical distributions of the re-scaled MLE at zero when sampling from the
Beta distribution (left) and the Gamma distribution (right): from top to bottom we have $\alpha=0.2, 0.5, 0.8$.   
}
\label{fig:cdfs}
\end{figure}

\section{Some Numerical Results}
\label{sec:Numerical}

Figure~\ref{fig:cdfs} gives plots of the empirical distributions of
$m=10000$ Monte Carlo samples from the distributions of 
$\widehat{f}_n (0)/ (C_2 n^{\alpha}/(1-\alpha))^{1/(1-\alpha)})$
when $n=200$ and $n=500$, together with the limiting distribution function obtained
in (\ref{LimitDistrib}).
The true density $f_0$ on the right side in Figure~\ref{fig:cdfs}  is
\begin{eqnarray}
f_0(x) = \int_0^{\infty} \frac{1}{y} 1_{[0,y]} (x)  \frac{y^{c -1}}{\Gamma (c)} \exp(-y) dy ;
\label{GammMix}
\end{eqnarray}
For $c \in (0,1)$, this family satisfies (G2) with $\alpha = 1-c$ and
$C_2 = 1/(\alpha \Gamma (1-\alpha))$.  (Note that for $c=1$, $f_0 (x) \sim \log (1/x)$ as $x \searrow 0$.)

The true density $f_0$ on the left side  in Figure~\ref{fig:cdfs}  is
\begin{eqnarray}
f_0(x) =  \frac{1}{\mbox{Beta}(1-a,2)} x^{-a} (1-x) 1_{(0,1]} (x) ;
\label{Beta}
\end{eqnarray}
For $a \in [0,1)$, this family satisfies (G2) with $\alpha = a$ and
$C_2 = 1/\mbox{Beta}(1-\alpha,2)$.

\begin{figure}[htb!]
\centering
\includegraphics[width=0.49\textwidth, height=0.37\textwidth]{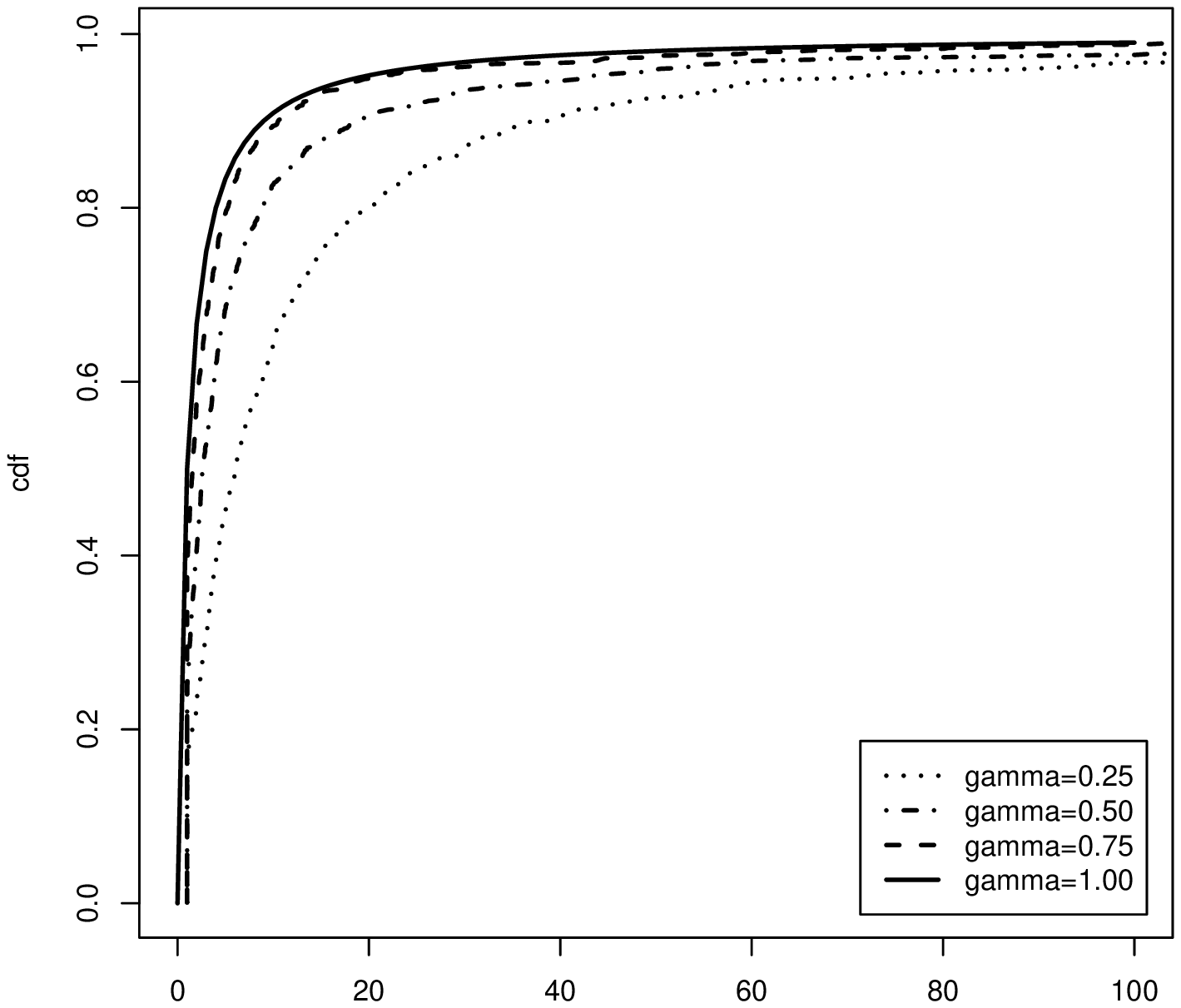}
\includegraphics[width=0.49\textwidth, height=0.37\textwidth]{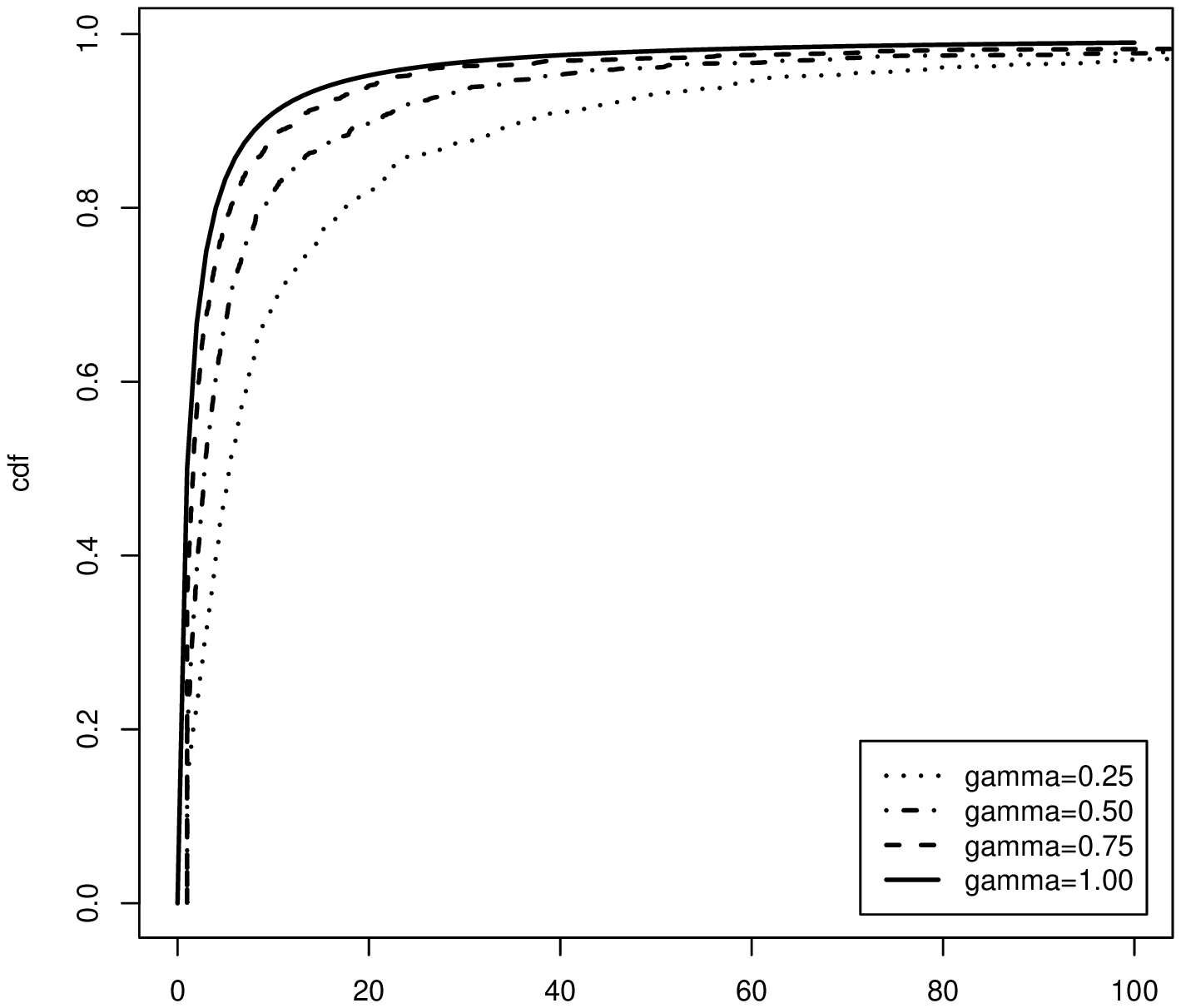}
\includegraphics[width=0.49\textwidth, height=0.37\textwidth]{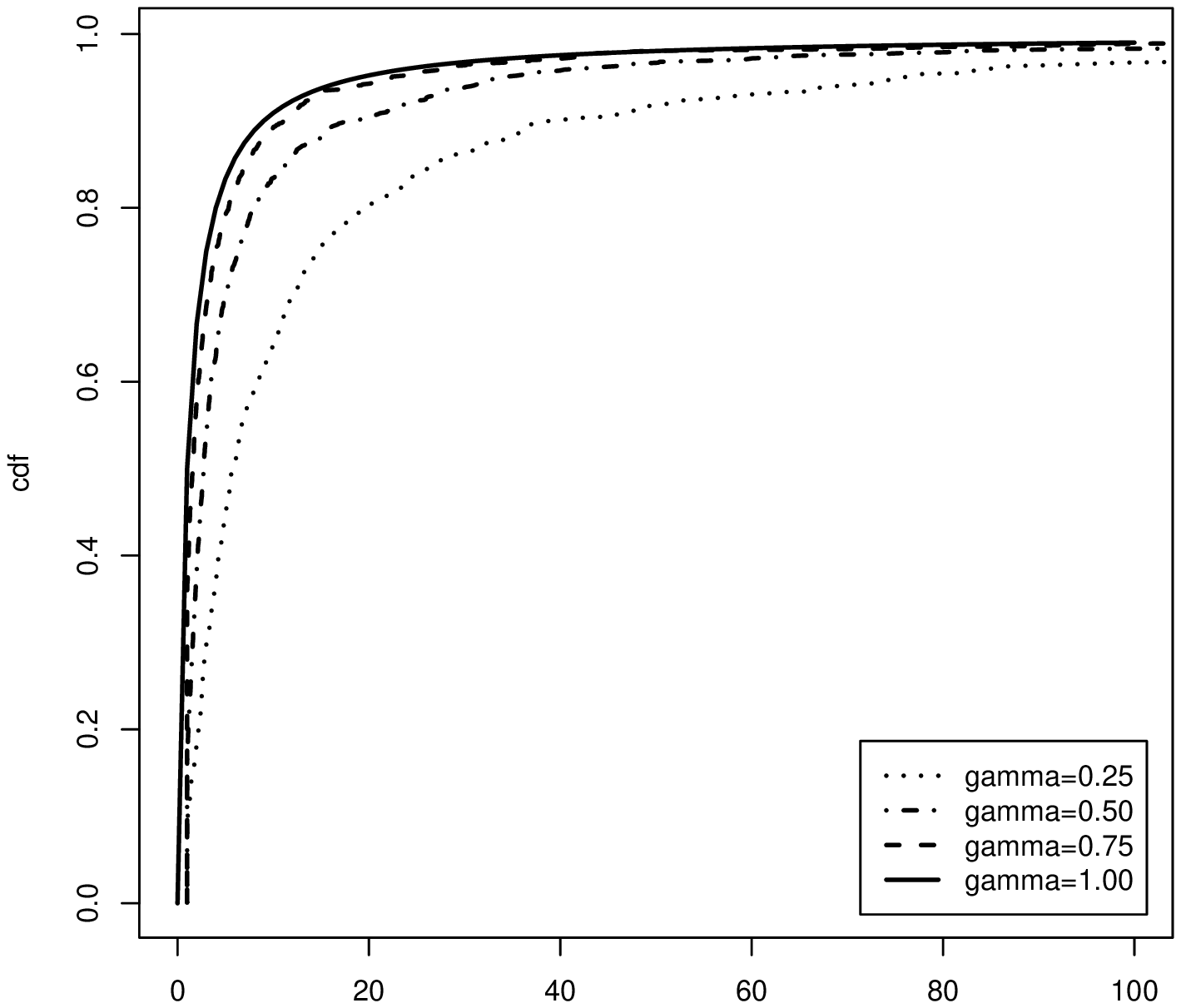}
\includegraphics[width=0.49\textwidth, height=0.37\textwidth]{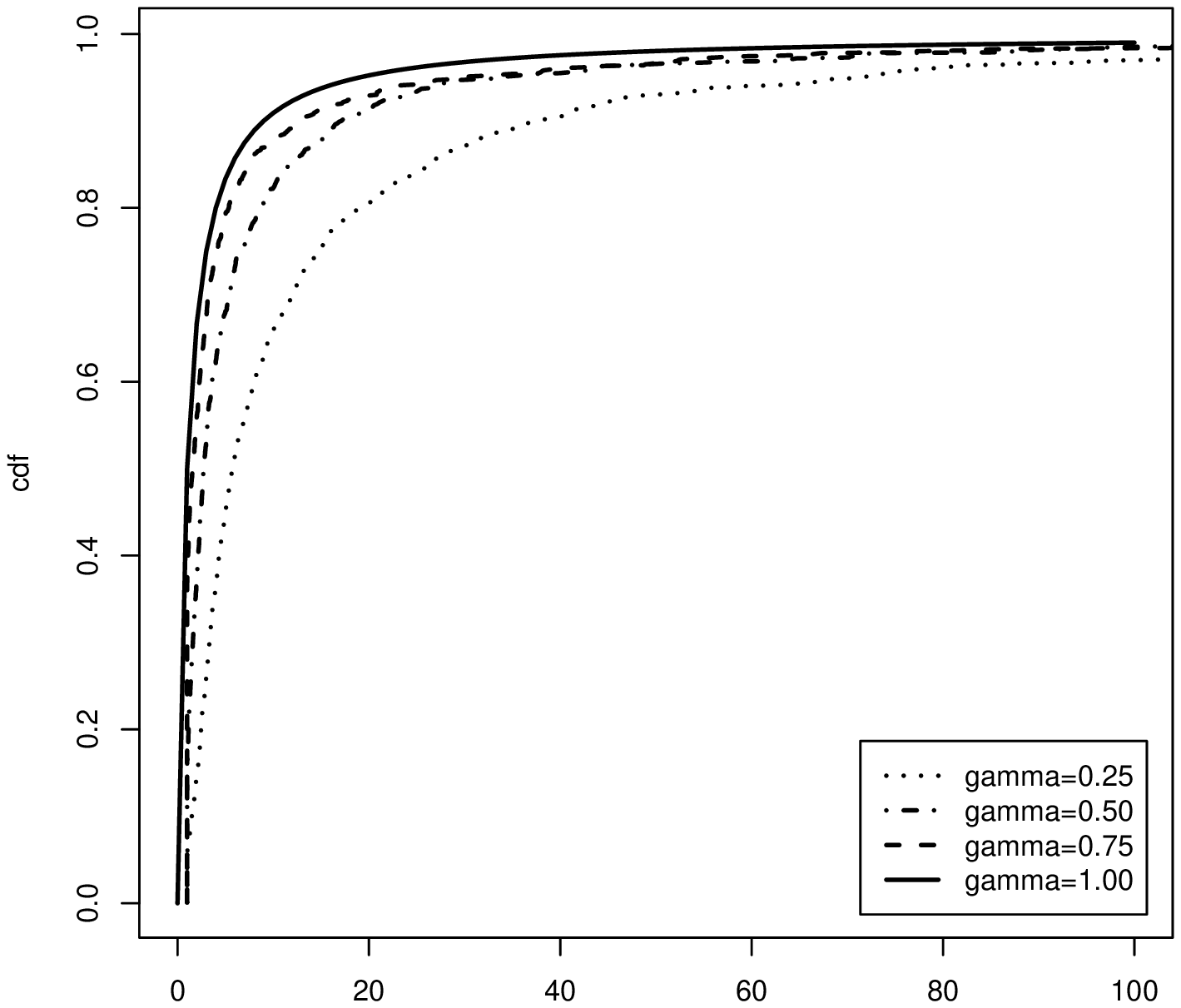}
\caption{Empirical distributions of the supremum measure: the cutoff values shown 
are $c=5$ (top left), $c=25$ (top right), $c=100$ (bottom left), $c=1000$ (bottom right).}
\label{fig:cdfsForCproblem}
\end{figure}

Figure~\ref{fig:cdfsForCproblem} shows simulations of the limiting distribution
\begin{eqnarray}\label{line:limquant}
\sup_{0\leq t \leq c}\left|t^{1-\gamma}\widehat h(t)/\gamma-1\right|
\end{eqnarray}
for different values of $c$ and $\gamma$.  Recall that if $\gamma=1$  the supremum
occurs at $t=0$ regardless of the value of $c$, and the limiting distribution
\eqref{line:limquant} has cumulative distribution function $1-1/(x+1).$
However, for $\gamma<1$, the distribution of \eqref{line:limquant}
depends both on $\gamma$ and on $c$, although the dependence
on $c$ is not visually prominent in Figure \ref{fig:cdfsForCproblem}.
Table \ref{tab:probone} shows estimated values of
\begin{eqnarray}\label{line:probone}
P\left(\sup_{0\leq t\leq c}|t^{1-\gamma}\widehat h(t)/\gamma-1|=1\right)
\end{eqnarray}
for different $c$ and $\gamma<1$, which clearly depends on the cutoff value $c$
(upper bound on the standard deviation in each case is 0.016).  Note that
\eqref{line:limquant} is equal to one if the location of the supremum occurs at
$t=0$ (with probability one).

\begin{table*}[htb!]
\caption{Simulation of \eqref{line:probone} for different values of $\gamma$ and $c$.}
\label{tab:probone}
\begin{tabular}{|c|ccccc|}
        \hline
         & $c=0.5$ & $c=5$ & $c=25$ & $c=100$ & $c=1000$\\
        \hline
        $\gamma=0.25$ & 0.361 & 0.171 & 0.140 & 0.092 & 0.06\\
        $\gamma=0.50$ & 0.422 & 0.249 & 0.190 & 0.162 & 0.148\\
        $\gamma=0.75$ & 0.489 & 0.387 & 0.349 & 0.358 & 0.367\\
        \hline
    \end{tabular}
\end{table*}
  
\begin{figure}[htb!]
\centering
\includegraphics[width=0.325\textwidth, height=0.3\textwidth]{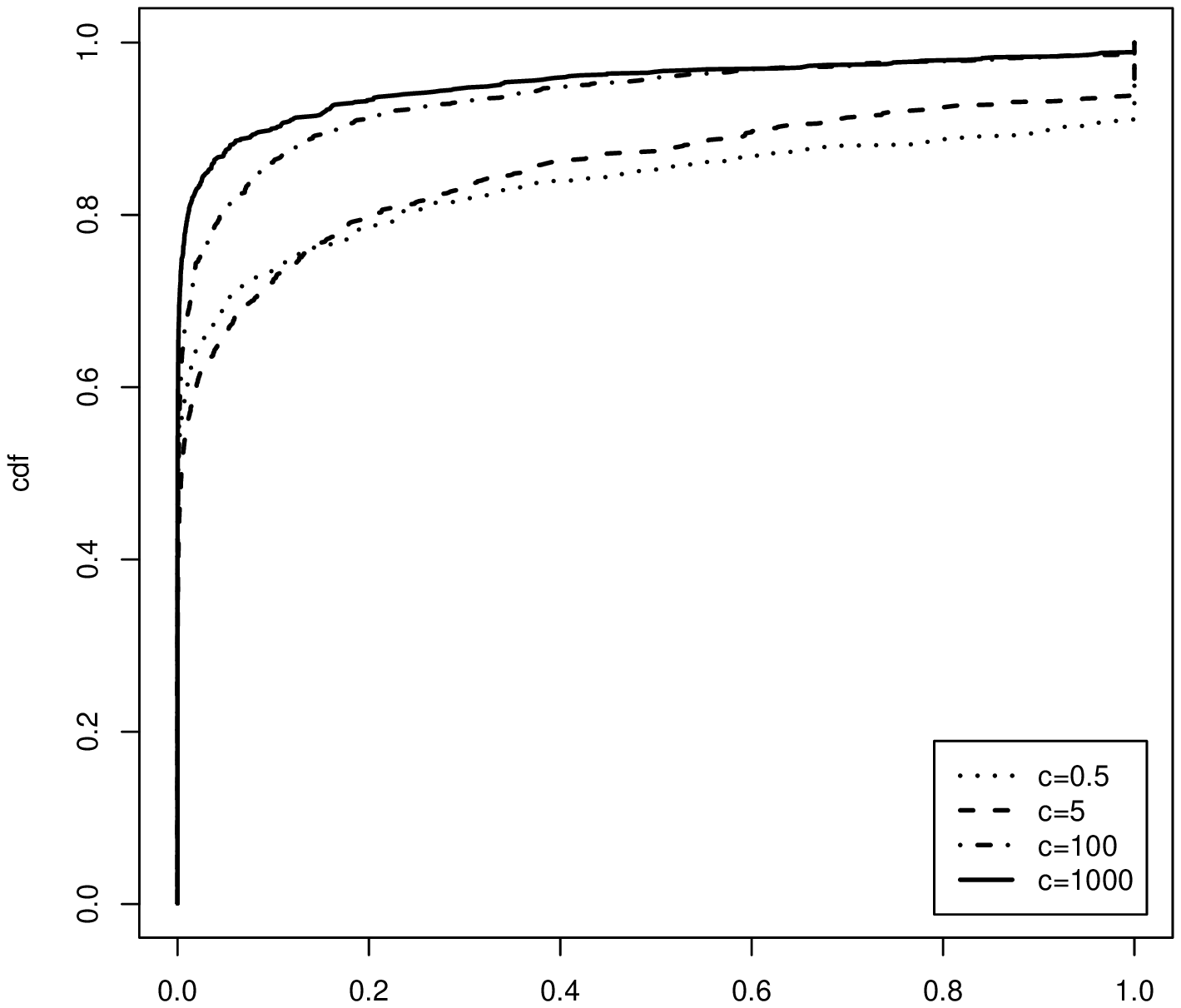}
\includegraphics[width=0.325\textwidth, height=0.3\textwidth]{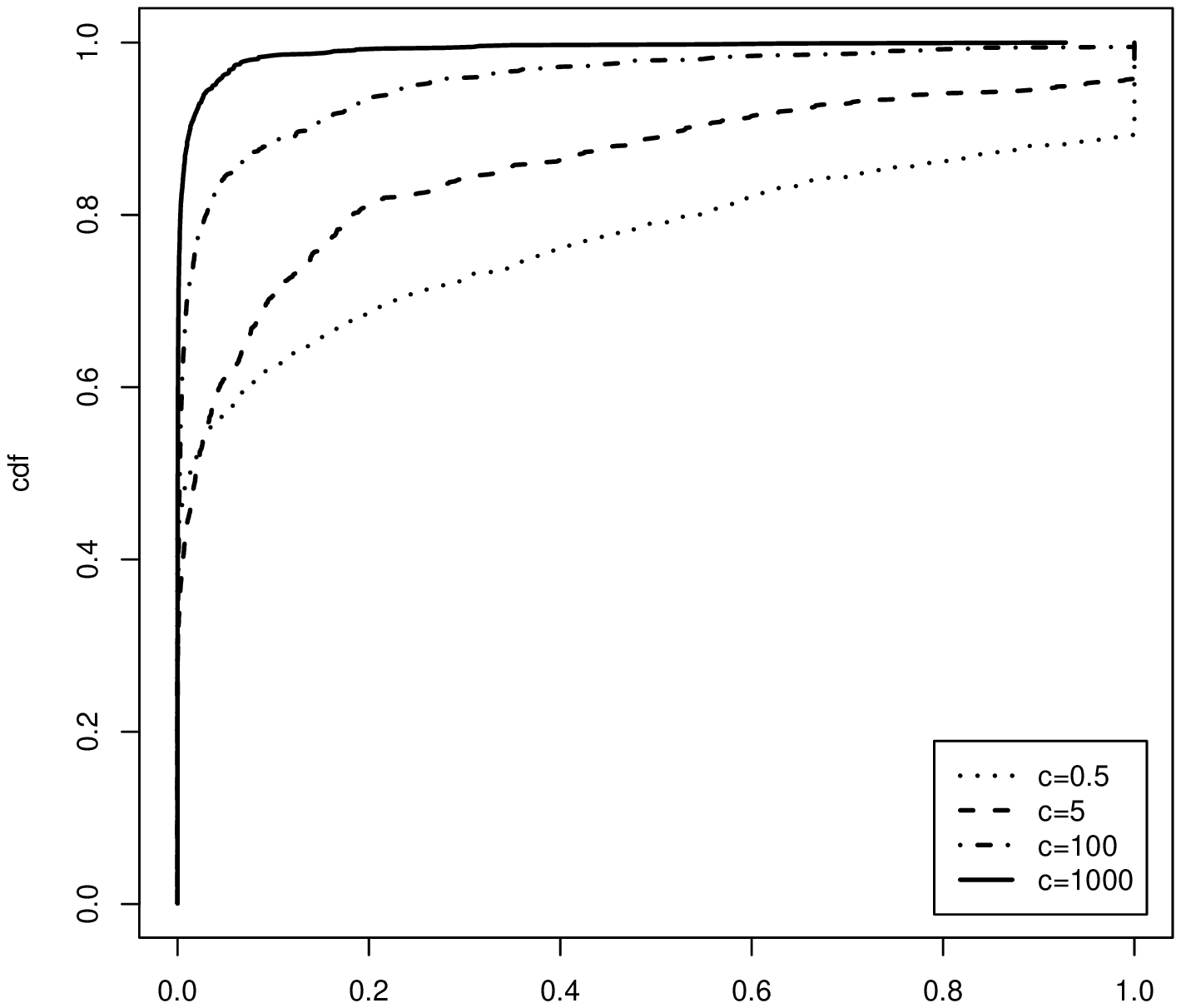}
\includegraphics[width=0.325\textwidth, height=0.3\textwidth]{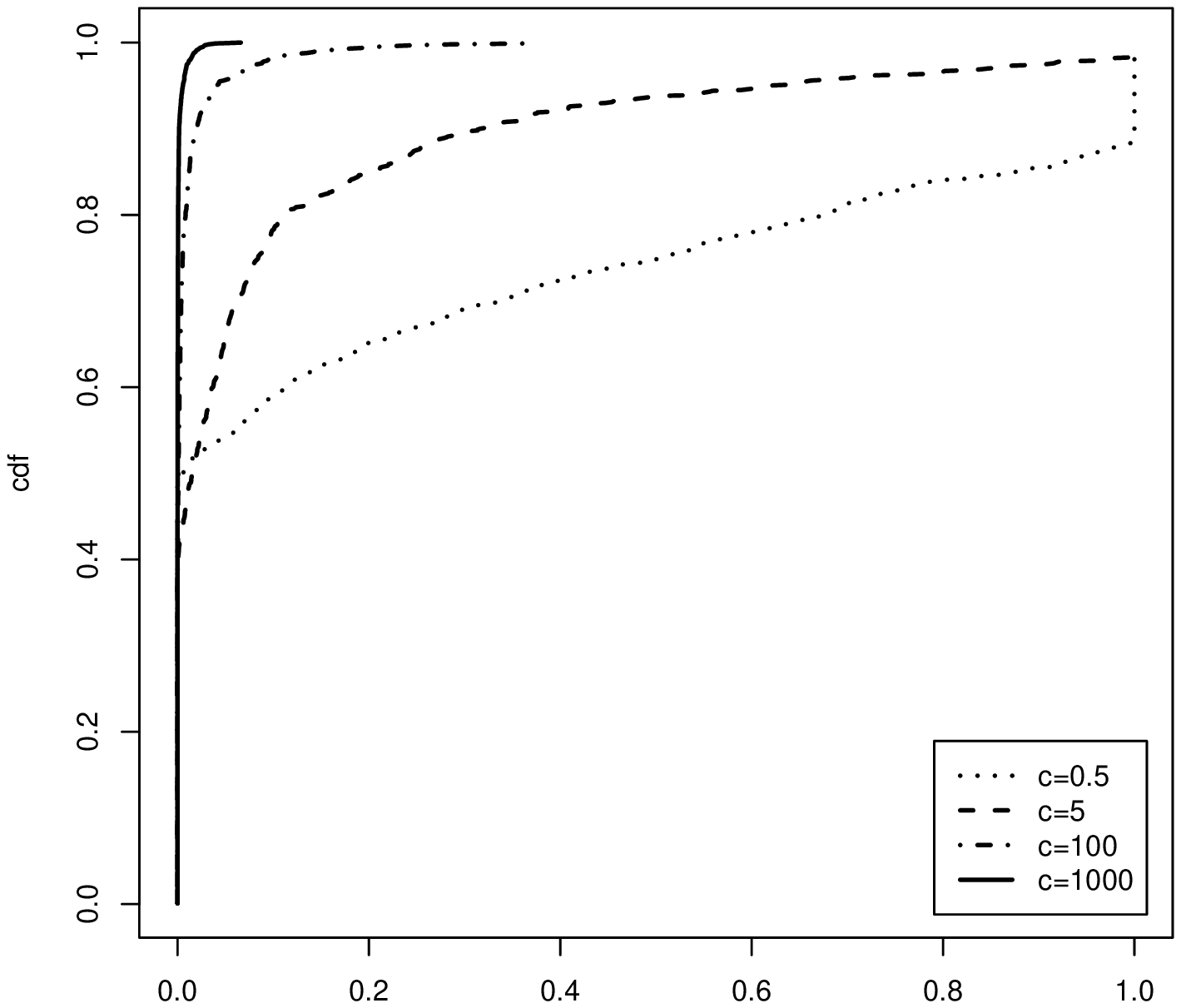}
\caption{Empirical distributions of the location where the supremum occurs:
from left to right we have $\gamma=0.25, 0.50, 0.75$.
Recall that for $\gamma=1$, the (non-unique) location of the supremum is
always zero by Corollary~\ref{cor:BoundedCase}.
The data were re-scaled to lie within the interval $[0,1]$.}
\label{fig:cdfs_loc}
\end{figure}

Cumulative distribution functions for the location of
the supremum in \eqref{line:limquant} are shown in Figure \ref{fig:cdfs_loc},
which clearly depend both on $\gamma$ and on $c$.

\section{Application to Mixtures}
\label{sec:MixtureApplications}

\subsection{Behavior near zero}
\label{BehaviorNearZero}

\begin{figure}[htb!]
\centering
\includegraphics[width=0.8\textwidth, height=0.5\textwidth]{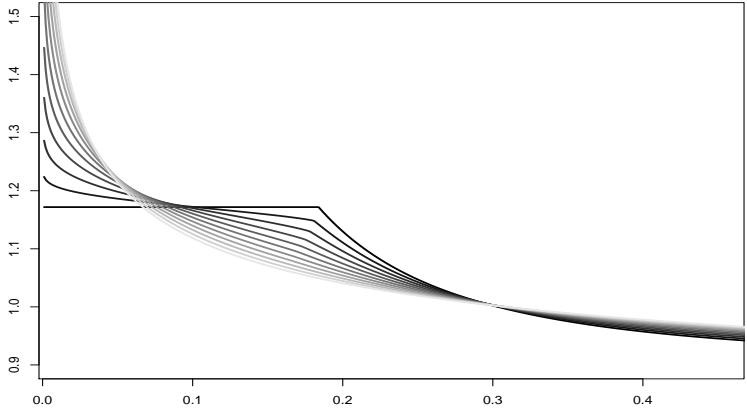}
\caption{Generalized Gaussian (or Subbotin) mixture densities with
$\epsilon = .1$, $\mu = 1$, $r \in \{ 1.0, 1.2, \ldots , 2.0\}$
  (black to light grey, respectively) as given by
(\ref{SubbotinMixtureModelTransformedToZeroOne}).}
\label{fig:GeneralizeNormalMixtureDensitiesTransformedToZeroOne}
\end{figure}

First, suppose that
$X_1, \ldots , X_n$ are i.i.d. with distribution function $F$ where,
\begin{eqnarray*}
&& \mbox{under} \ H_0 : \ F = \Phi_r, \qquad \mbox{the generalized normal distribution} \\
&& \mbox{under} \ H_1 : \ F = (1-\epsilon) \Phi_r+ \epsilon \Phi_r (\cdot - \mu ) ,  \ \ \epsilon \in (0,1),  \ \ \mu > 0,
\end{eqnarray*}
where $\Phi_r (x) \equiv \int_{-\infty}^x \phi_r (y) dy$ with $\phi_r (y) \equiv \exp ( - |y|^r / r )/C_r$ for $r>0$ gives
the generalized normal (or Subbotin) distribution;
here $C_r \equiv 2 \Gamma (1/r) r^{(1/r)-1}$ is the normalizing constant.
If we transform to $Y_i \equiv 1-\Phi_r (X_i) \sim G$, then, for $0~\le~y~\le~1$,
\begin{eqnarray*}
&& \mbox{under} \ H_0 : G(y) = y,  \qquad \mbox{the Uniform}(0,1)\ \mbox{d.f.}, \\
&& \mbox{under} \ H_1 :  G(y) = G_{\epsilon, \mu, r} (y) = (1-\epsilon) y + \epsilon (1- \Phi_r ( \Phi_r^{-1}(1-y) - \mu )).
\end{eqnarray*}
Let $g_{\epsilon, \mu,r}$ denote the density of  $G_{\epsilon, \mu,r}$; thus 
\begin{eqnarray}
\hspace{1cm}g_{\epsilon, \mu,r} (y)
& = & 1-\epsilon + \epsilon \exp \left \{ - \frac{1}{r} \left ( | \Phi_r^{-1} (1-y) - \mu |^r - | \Phi_r^{-1} (1-y) |^r \right ) \right \}
          \label{SubbotinMixtureModelTransformedToZeroOne} .
\end{eqnarray}
It is easily seen that $g_{\epsilon, \mu,r}$ is monotone decreasing on $(0,1)$ and is unbounded at zero if $r>1$.
Figure~\ref{fig:GeneralizeNormalMixtureDensitiesTransformedToZeroOne}
shows plots of these
densities for $\epsilon = .1$, $\mu = 1$, and $r \in \{ 1.0, 1.1, \ldots , 2.0 \}$.
Note that $g_{\epsilon, \mu,1}$ is
bounded at $0$:  in fact $g_{\epsilon, \mu, 1} (y) = 1-\epsilon + \epsilon e^{\mu}$ for $0 \le y \le 2^{-1} e^{- \mu}$.

\begin{prop}
The distribution $F_{\mu,r} (y) \equiv 1 - \Phi_r (\Phi_r^{-1} (1-y) - \mu ) $ is regularly varying at $0$ with exponent $1$.
That is, for any $c>0$,
$$
\lim_{y \rightarrow 0+} \frac{F_{\mu, r} (cy)}{F_{\mu, r} (y)} = c,
$$
i.e. (\ref{GnedenkoHypothesis}) holds 
with $\gamma = 1$.
\end{prop}
\medskip

\par\noindent
{\bf Proof.} 
Define $\kappa_r(y)=\Phi_r^{-1}(1-y).$  Our first goal will be to show that 
\begin{eqnarray}\label{line:kappa1}
\lim_{y\rightarrow0}\frac{\kappa_r(y)}{\tilde{\kappa}_r (y)}=1,
\end{eqnarray}
where (for $y$ small)
\begin{eqnarray*}
\tilde{\kappa}_r (y) =
\left (-r \log \left ({C_r \ y \ \left\{ r\log \left(\frac{1}{C_r y} \right)\right\}^{(r-1)/r} } \right ) \right )^{1/r}.
\end{eqnarray*}

To prove \eqref{line:kappa1}, it is enough to show that 
\begin{eqnarray}\label{line:kappa2}
\lim_{y\rightarrow0}{\tilde{\kappa}_r (y)^{r-1}}(\kappa_r(y)-\tilde{\kappa}_r (y))=0.
\end{eqnarray}
This result follows from  \citet[Theorem 1.1.2]{MR2234156}.  Define
\begin{eqnarray*}
b_n = \tilde\kappa_r(1/n), \ \ \ a_n = 1/b_n^{r-1},
\end{eqnarray*}
and choose $F=\Phi_r$ in the statement of Theorem 1.1.2.  Then, if we can show that 
\begin{eqnarray}\label{ExtremeValueLimitTailCond}
n (1 - \Phi_r ( a_n x + b_n )) \rightarrow \log G(x)\equiv e^{-x}  ,  \qquad x \in \RR,
\end{eqnarray}
it follows from \citet[Theorem 1.1.2 and Section 1.1.2]{MR2234156} that  for all $x\in\RR$
\begin{eqnarray*}
\lim_{y\rightarrow 0} \frac{U(x/y)-b_{\lfloor1/y\rfloor}}{a_{\lfloor1/y\rfloor}} = G^{-1}(e^{-1/x})= \log(1/x),
\end{eqnarray*}
where $U(t)=(1/(1-\Phi_r))^{-1}(t)= \Phi_r^{-1}(1-1/t).$  Choosing $x=1$ yields \eqref{line:kappa2}.  Therefore, we need to prove \eqref{ExtremeValueLimitTailCond}.

To do this, we make use of the following, which is a generalization of Mills' ratio to the generalized Gaussian family
\begin{eqnarray}
1-\Phi_r (z) \sim \frac{\phi_r (z)}{z^{r-1}}  \qquad \mbox{as} \ \ z \rightarrow \infty .
\label{GeneralizationOfMillsRatio}
\end{eqnarray}
The statement follows from l'H\^{o}pital's rule:
\begin{eqnarray*}
\lim_{z \rightarrow \infty} \frac{\int_z^\infty \phi_r (y) dy}{z^{1-r} \phi_r (z)}
& = & \lim_{z \rightarrow \infty} \frac{-\phi_r (z)}{(1-r)z^{-r} \phi_r (z) + z^{1-r} \phi_r (z) (-z^{r-1} )} \\
& = & \lim_{z \rightarrow \infty} \frac{1}{1 - (1-r) z^{-r} } = 1 .
\end{eqnarray*}
Now, 
\begin{eqnarray*}
n(1- \Phi_r (a_n x + b_n ) )
& \sim & n \frac{\phi_r (a_n x + b_n )}{(a_n x + b_n )^{r-1} } \\
& = & \frac{n}{C_r b_n^{r-1}}
          \frac{\exp\left ( - \frac{b_n^r}{r} \left ( 1 + \frac{a_n x}{b_n} \right )^r \right )}{(1+ a_n x/b_n )^{r-1}} \\
& \sim & \frac{n}{C_r b_n^{r-1}} \exp\left ( - \frac{b_n^r}{r} \left (1 + \frac{r  x}{b_n^r} \right ) \right ) \\
& = & \exp \left ( - \left ( \frac{b_n^r}{r} + (r-1) \log b_n - \log n + \log C_r \right ) \right ) \exp ( - x)  \\
& \rightarrow & \exp (- 0) \cdot \exp (-x)
\end{eqnarray*}
by using the definition of $b_n$.  We have thus shown that \eqref{line:kappa1} holds. 
\medskip

Then, for $y\rightarrow 0$, by \eqref{GeneralizationOfMillsRatio} and \eqref{line:kappa1} 
\begin{eqnarray*}
F_{\mu, r}(y)& = & 1- \Phi_r ( \kappa_r (y) - \mu ) \sim 1- \Phi_r ( \tilde\kappa_r (y) - \mu ) \\
&\sim& \frac{\phi_r(\tilde{\kappa}_r (y)-\mu)}{(\tilde{\kappa}_r (y)-\mu)^{r-1}}.
\end{eqnarray*}
Plugging in the definition of $\phi_r,$ we find that
\begin{eqnarray*}
F_{\mu, r}(y)& \sim & \frac{1/C_r}{(\tilde{\kappa}_r (y)-\mu)^{r-1}}
            \exp\left(-\frac{\tilde{\kappa}_r (y)^r}{r}\left|1-\frac{\mu}{\tilde{\kappa}_r (y)}\right|^r\right)\\
&=& \frac{1/C_r}{(\tilde{\kappa}_r (y)-\mu)^{r-1}}
        \exp\left\{ \left ( \log (C_r y)  + \log(r \log (1/(C_r y)) ) \right ) \left|1-\frac{\mu}{\tilde{\kappa}_r (y)}\right|^r\right)\\
&=& \frac{1/C_r}{(\tilde{\kappa}_r (y)-\mu)^{r-1}}\ (C_r y)^{\left|1-\frac{\mu}{\tilde{\kappa}_r (y)}\right|^r}
         \cdot \left \{ r \log \frac{1}{C_r y} \right \}^{\frac{r-1}{r} \left|1-\frac{\mu}{\tilde{\kappa}_r (y)}\right|^r } .
\end{eqnarray*}
Note that $\lim_{y\rightarrow 0} \tilde\kappa_r(cy)/\tilde\kappa_r(y)=1$.  Therefore,
\begin{eqnarray*}
\frac{F_{\mu, r}(cy)}{F_{\mu, r}(y)}
&\sim& c^{\left|1-\frac{\mu}{\tilde{\kappa}_r (cy)}\right |^r}
         \cdot (C_r y)^{{\left|1-\frac{\mu}{\tilde{\kappa}_r (cy)}\right|^r}-{\left|1-\frac{\mu}{\tilde{\kappa}_r (y)}\right|^r}}
         \cdot \left(\frac{\tilde{\kappa}_r (y)-\mu}{\tilde{\kappa}_r (cy)-\mu}\right)^{r-1} \\
&& \qquad \cdot  \frac{\left \{ r \log \frac{1}{C_r c y} \right \}^{\frac{r-1}{r} \left|1-\frac{\mu}{\tilde{\kappa}_r (cy)}\right|^r }}
                                  {\left \{ r \log \frac{1}{C_r y} \right \}^{\frac{r-1}{r} \left|1-\frac{\mu}{\tilde{\kappa}_r (y)}\right|^r }}\\
&\rightarrow & c \cdot 1 \cdot 1 \cdot 1 = c .
\end{eqnarray*}
Thus (\ref{GnedenkoHypothesis}) holds with $\gamma = 1$.
\hfill $\Box$
\medskip

By the theory of regular variation (see e.g. \citet[page 21]{MR1015093}), this implies that
$F_{\mu, r} (y) = y \ell (y)$ where $\ell$ is slowly varying at $0$.  It then follows easily that
(\ref{GnedenkoHypothesis}) holds for $F_0 = G_{\epsilon, \mu, r}$ with exponent $1$.
Thus our theory of Section 1 applies with $a_n$ of Theorem 1.1 taken to be $a_n = G_{\epsilon, \mu , \gamma} (1/n)$;
i.e.
$$
\frac{1}{n} = G_{\epsilon, \mu, r} (a_n ) = (1-\epsilon )a_n  + \epsilon F_{\mu , r} (a_n )
\ \dot= \ \epsilon F_{\mu, r} (a_n )
$$
where the last approximation is valid for $r>1$, but not for $r=1$.   When $r=1$, the
first equality can be solved explicitly, and we find:
\begin{eqnarray}
a_n = \left \{
\begin{array}{l l} 1 - \Phi_r ( \Phi_r^{-1} (1 - (1/(n\epsilon))) + \mu ) , & \mbox{when}\ \ r> 1\\
                           n^{-1} (1-\epsilon + \epsilon e^{\mu} )^{-1} , & \mbox{when} \ \ r = 1.
\end{array} \right .
\label{NormalizingSequenceSubbotinDistributions}
\end{eqnarray}
We conclude that Theorem~\ref{thm:GeneralTheoremGrenanderAtZero} holds for $a_n $ as in the last display
where $\widehat{f}_n$ is the Grenander estimator of $g_{\epsilon, \mu,r}$ based on  $Y_1 , \ldots , Y_n$.

Another interesting mixture family to consider is as follows:
suppose that $\Phi_1$, $\Phi_2$ are two fixed distribution functions:  then
\begin{eqnarray*}
&& \mbox{under} \ H_0 : \ F = \Phi_1,  \\
&& \mbox{under} \ H_1 : \ F = (1-\epsilon) \Phi_1+ \epsilon \Phi_2  ,  \ \ \epsilon \in (0,1) .
\end{eqnarray*}
Using the transformation to $Y_i \equiv 1-\Phi_1 (X_i) \sim G$, then, for $0 \le y \le 1$ we find that under
$H_1$ the distribution of the $Y_i$'s is given by
\begin{eqnarray*}
&& G(y) = (1-\epsilon) y  + \epsilon (1 - \Phi_2 (\Phi_1^{-1} (1-y))), \\
&& g(y) = (1-\epsilon) + \epsilon \frac{\phi_2 ( \Phi_1^{-1} (1-y))}{\phi_1 (\Phi_1^{-1} (1-y))} .
\end{eqnarray*}
For $\Phi_2 $ given in terms of $\Phi_1$ by the (Lehmann alternative) distribution function
$\Phi_2 (y) = 1 - (1- \Phi_1 (y))^{\gamma}$, this becomes
\begin{eqnarray*}
&& G(y) = (1-\epsilon) y  + \epsilon y^{\gamma}, \\
&& g(y) = (1-\epsilon) + \epsilon \gamma y^{\gamma-1}  .
\end{eqnarray*}
When $0 < \gamma < 1$
 this family fits into the framework of our condition {\bf G2} with $\alpha = 1- \gamma$ and $C_2 = \epsilon \gamma$.

\subsection{Estimation of the contaminating density}
\label{ContaminationDensityEstimation}

Suppose that
$G_{\epsilon, F} (y) = (1-\epsilon) y + \epsilon F(y)$
where $F$ is a concave distribution on $[0,1]$ with monotone decreasing density $f$.
Thus the density $g_{\epsilon, F}$  of $G_{\epsilon , F}$ is given by
$g_{\epsilon, F} (y) = (1-\epsilon ) + \epsilon f(y)$.
Note that $g_{\epsilon , F}$ is also monotone decreasing,
and $g_{\epsilon, F} (y) \ge 1-\epsilon + \epsilon f(1) = 1-\epsilon = g_{\epsilon,F} (1)$ if $f(1) = 0$.
For $\epsilon>0$ we can write
$$
f(y) = \frac{g_{\epsilon,F} (y) - (1-\epsilon)}{\epsilon} .
$$
If $Y_1 , \ldots , Y_n$ are i.i.d. $g_{\epsilon, F}$ then we can estimate $g_{\epsilon, F}$ by the Grenander estimator
$\widehat{g}_n$, and we can estimate $\epsilon$ by
$$
\widehat{\epsilon}_n = 1- \widehat{g}_n (1) .
$$
 This results in
the following estimator $\widehat{f}_n$ of the contaminating density $f$:
\begin{eqnarray*}
\widehat{f}_n (y)
 = \frac{\widehat{g}_n (y) - (1-\widehat{\epsilon}_n )}{\widehat{\epsilon}_n }
  =  \frac{\widehat{g}_n (y) - \widehat{g}_n (1) }{1 - \widehat{g}_n (1) } ,
\end{eqnarray*}
which is quite similar in spirit to a setting studied by
\citet{MR1701099}.  Here, however, we propose using the shape constraint of monotonicity,
and hence the Grenander estimator, to estimate both $\epsilon$ and $f$.
We intend to study this estimator elsewhere.

\setcounter{section}{5}
\section*{Appendix A: Proofs for Section 1}
\label{sec:proofs}

Before proving Theorem~\ref{thm:GeneralTheoremGrenanderAtZero}, we need the following two
lemmas.  The first lemma shows that the functionals $\argmax^R$ and $\argmax^L$
are both $O_p (1)$, while the second shows these are equivalent almost surely for the limiting
Poisson process.
Together, these two lemmas will show that both functionals $\argmax^R$ and $\argmax^L$ are continuous.
Below we assume that (\ref{GnedenkoHypothesis}) holds and that $n F_0 (a_n ) \sim 1$.  
Thus both (\ref{SequentialConvergenceRightTailCond}) and (\ref{GeneralWeakConvergenceECPGnedenkoHypothesis}) 
also hold.

\begin{lem}
\label{lem:BigOhPForArgmax}
(i) \ When $\gamma =1$ and $x>1$, $\mbox{argmax}_v^{L,R} \{ n \FF_n (a_n v) - x v \} = O_p (1)$.\\
(ii) \ When $\gamma \in (0,1)$ and $x>0$, $\mbox{argmax}_v^{L,R} \{ n \FF_n (a_n v ) - xv \} = O_p (1)$.
\end{lem}

\par\noindent
{\bf Proof.}   
It suffices to show that
\begin{eqnarray*}
\limsup_{n\rightarrow \infty} P( \sup_{v \ge K} \{ n \FF_n (a_n v) - x v\} \ge 0) \rightarrow 0,
\ \ \mbox{as} \ \ K \rightarrow \infty
\end{eqnarray*}
under the conditions specified.
Let $h(x) = x(\log x-1)+1$ and recall the inequality
$$
P( \mbox{Bin}(n,p)/(np) \ge t) \le \exp ( - n p h(t))
$$
for $t\ge1$ where $\mbox{Bin}(n,p)$ denotes a Binomial$(n,p)$ random variable;
see e.g.
\citet[inequality 10.3.2, page 415]{MR838963}.
It follows that
\begin{eqnarray}
\lefteqn{P( \sup_{v \ge K} \{ n \FF_n (a_n v) - x v\} \ge 0) } \nonumber \\
& = & P( \cup_{j=K}^{\infty} \{ n \FF_n (a_n v) - x v \ge 0 \ \ \mbox{for some} \ \ v \in [j, j+1) \})
         \nonumber   \\
& \le & \sum_{j=K}^\infty P( n \FF_n ( a_n (j+1)) - x j \ge 0 ) \nonumber \\
& = & \sum_{j=K}^\infty P\left ( \frac{ n \FF_n (a_n (j+1))}{n F_0 (a_n (j+1) )} \ge \frac{x j }{n F_0 (a_n (j+1))} \right )
          \nonumber  \\
& \le & \sum_{j=K}^\infty \exp \left ( - n F_0 ( a_n (j+1)) h\left (\frac{xj}{nF_0 (a_n (j+1))} \right ) \right )
\label{line:lastbound}
\end{eqnarray}
Next, since $F_0$ is concave,
$$
n F_0(a_n(j+1))\leq nF_0(a_n(K+1)) \frac{j+1}{K+1}
$$
for $j\geq K$ and $nF_0(a_n(K+1))\rightarrow (K+1)^\gamma$ and $n\rightarrow \infty$.
Therefore, for all $j\geq K$ and sufficiently large $n$, we have
\begin{eqnarray*}
\frac{xj}{ nF_0(a_n(j+1))} \geq \delta (K+1)^{1-\gamma} \frac{xj}{j+1}
\end{eqnarray*}
for any fixed $\delta<1.$  We need to handle the two cases $\gamma=1$ and
$\gamma<1$ separately.  Note that if $\gamma<1$, then the above display
shows that $K,n$ can be chosen sufficiently large so that $(xj)/nF_0(a_n(j+1))$
is uniformly large.  On the other hand, if $\gamma =1$ and $x>1$ then we can
pick $\delta, K,n$ large enough so that $(xj)/nF_0(a_n(j+1))$ is strictly greater
than $1+\eps$ for some $\eps>0$, again uniformly in $j$.

Suppose first that $\gamma<1$.  Then for $K,n$ large, since $h(x)\sim x \log x $
as $x \rightarrow \infty$, there exists a constant $0<C<1$ such that for all $j \ge K$
\begin{eqnarray*}
n F_0 ( a_n (j+1)) h\left (\frac{xj}{nF_0 (a_n (j+1))} \right )
&\geq & C (xj)\log\left(\frac{xj}{j+1}\right)\\ 
&\geq& C_{x}(xj),
\end{eqnarray*}
for some other constant $C_{x}>0$.  This shows that the sum in
(\ref{line:lastbound}) converges to zero as $K\rightarrow \infty,$ as required.

Suppose next that $\gamma=1$.  Note that the function $h(x)>0$  for $x>1$.
Therefore, combining our arguments above, we find that for all $j\ge K$
\begin{eqnarray*}
n F_0 ( a_n (j+1)) h\left (\frac{xj}{nF_0 (a_n (j+1))} \right )
&\geq & \delta (j+1) h\left (\frac{xj}{nF_0 (a_n (j+1))} \right )\\
&\geq & C_{x, \delta}(j+1),
\end{eqnarray*}
again for some $C_{x,\delta}>0$.  This again implies that the sum in
(\ref{line:lastbound}) converges to zero as $K\rightarrow \infty,$ and completes the proof.
\hfill $\Box$

\begin{lem}
\label{lem:ArgmaxUniqueForPoissonLimit}
Suppose that $\gamma \in (0,1]$.  Then
$$
V_x^{L} \equiv \argmax_v^L \{ \NN (v^{\gamma} ) - x v \} 
= \argmax_v^R \{ \NN (v^{\gamma} ) - x v \}  \equiv V_x^R
\qquad a.s.
$$
\end{lem}

\par\noindent
{\bf Proof.}  Suppose that $V_x^L < V_x^R$.  Then it follows that
$\NN ( (V_x^L)^{\gamma} ) - x V_x^L = \NN ( (V_x^R)^{\gamma} ) - x V_x^R$,
or, equivalently
$$
\NN( (V_x^R)^{\gamma} ) - \NN ( (V_x^L)^{\gamma} ) = x \{ V_x^R - V_x^L \} .
$$
Now $(V_x^R)^{\gamma} , \ (V_x^L)^{\gamma} \in J ( \NN) \equiv \{ t>0 : \ \NN(t) - \NN(t-) \ge 1 \}$,
so the left side of the last display takes values in the set $ \{ 1, 2 , \ldots \}$, while the right side
takes values in $x \cdot \{ r^{1/\gamma} - s^{1/\gamma} :  \ r , s \in J(\NN),  r> s \}$.  But it is well-known
that all the (joint) distributions of the points in $J(\NN) $ are absolutely continuous with respect to
Lebesgue measure, and hence the equality in the last display holds only for sets with probability $0$.
\hfill $\Box$
\smallskip

\par\noindent
{\bf Proof of Theorem~\ref{thm:GeneralTheoremGrenanderAtZero}}:
We first prove convergence of the one-dimensional distributions
of $ na_n \widehat{f}_n ( a_n t) $.
Fix $K>0$, and let $x > 1_{\{\gamma=1\}}$ and $t \in (0,K]$.
By the switching relation (\ref{ClosedSwitchingRelationGrenander}),
\begin{eqnarray*}
P( na_n \widehat{f}_n ( a_n t ) \le x )
& = & P( \widehat{s}_n^L(x/(na_n )) \le a_n  t ) \\
& = & P(  \mbox{argmax}_s^L \{ \FF_n (s) - xs/(na_n) \} \le a_n t ) \\
& = & P( \mbox{argmax}_v^L \{ \FF_n (v a_n) - x (v/n) \} \le t ) \\
& = & P( \mbox{argmax}_v^L \{ n \FF_n (v a_n ) - x v \} \le t  ) \\
& \rightarrow & P( \mbox{argmax}_v^L \{ \NN( v^{\gamma}) - xv \} \le  t ) \\
& = & P( \widehat{h}_{\gamma} (t) \le x )
\end{eqnarray*}
where the convergence
follows from (\ref{GeneralWeakConvergenceECPGnedenkoHypothesis}),
and the argmax continuous mapping theorem for $D[0,\infty)$ applied to the processes
$\{ v \mapsto n \FF_n (v a_n) - x v :  \ v \ge 0 \}$; see e.g.
\citet[Theorem 3 and Corollary 1]{MR2042258}.   Note that
Lemma~\ref{lem:BigOhPForArgmax} yields the
$O_p(1)$ hypothesis of Ferger's Corollary 1, while
Lemma~\ref{lem:ArgmaxUniqueForPoissonLimit} shows that equality
holds in the limit conclusion.

Convergence of the finite-dimensional distributions of
$ \widehat{h}_n  (t) \equiv na_n \widehat{f}_n ( a_n t)$
follows in the same way by using the process convergence in
(\ref{GeneralWeakConvergenceECPGnedenkoHypothesis})
for finitely many values $(t_1, x_1), \ldots , (t_m , x_m)$
where each $t_j \in \RR^+$ and $x_j > 1_{\{ \gamma = 1 \}}$.

 To verify tightness of $\widehat{h}_n $ in $D[0,\infty)$ we use
\citet[Theorem 16.8]{MR1700749}.  Thus, it is sufficient to show that for any $K>0$, and any $\eps>0$
\begin{eqnarray}
\lim_{M\rightarrow \infty}\limsup_n P\left(\sup_{0\leq t \leq K}|\widehat h_n(t)|\geq M\right)&=&0 \label{line:tight_condA}\\
\lim_{\delta\rightarrow 0}\limsup_n P\left( w_{\delta,K}(\widehat h_n)\geq \eps\right)&=&0 \label{line:tight_condB},
\end{eqnarray}
where $w_{\delta,K}(h)$ is the modulus of continuity in the Skorohod topology defined as
\begin{eqnarray*}
w_{\delta,K}(h) &=& \inf_{\{t_i\}_r} \max_{0<i\leq r} \sup\left\{|h(t)-h(s)|: s,t \in [t_{i-1}, t_i)\cap[0,K]\right\},
\end{eqnarray*}
where $\{t_i\}_r$ is a partition of $[0,K]$ such that $0=t_0 < t_1 < \ldots < t_r =K$
and $t_i-t_{i-1}>\delta$.   Suppose then that $h$ is a piecewise constant
function with discontinuities occurring at the (ordered) points $\{\tau_i\}_{i\geq 0}$.
Then if $\delta \leq \inf_i |\tau_i-\tau_{i-1}|$ we necessarily have that $w_{\delta,K}(h) =0.$

First, note that
since $\widehat{h}_n $ is non-increasing,
$$
\| \widehat{h}_n \|_0^m \equiv \sup _{0 \le t \le m} | \widehat{h}_n (t) | = \widehat{h}_n (0),
$$
and hence \eqref{line:tight_condA} follows from the finite-dimensional convergence proved above.

Next, fix $\eps>0$.  Let $0 =\tau_{n,0}< \tau_{n,1} < \cdots < \tau_{n,K_n}< K$
denote the (ordered) jump points of $\widehat{h}_n $, and let
$0= T_{n,0}< T_{n,1} < \cdots < T_{n,J_n}< K$ denote the (again, ordered)
jump points of $n\FF_n(a_n t)$.  Because
$ \{ \tau_{n,1} , \ldots , \tau_{n,K_n}\} \subset \{ T_{n,1}  , \ldots  , T_{n,J_n} \} $,
it follows that $\inf\{\tau_{i,n}-\tau_{i-1,n}\}\geq \inf\{T_{i,n}-T_{i-1,n}\}$ and hence
\begin{eqnarray*}
P\left( w_{\delta,K}(\widehat h_n)\geq \eps\right) \leq P\left( \inf_{i=1, \ldots, J_n}\{T_{i,n}-T_{i-1,n}\} < \delta\right).
\end{eqnarray*}
Now, by (\ref{GeneralWeakConvergenceECPGnedenkoHypothesis}) 
and continuity of the inverse map (see e.g. \citet[Theorem 13.6.3, page 446]{MR1876437})
\begin{eqnarray*}
(T_{n,1} , \ldots , T_{n,J_n} , 0 , 0 , \ldots ) \Rightarrow ( T_1^{1/\gamma}, \ldots , T_J^{1/\gamma} , 0 , 0, \ldots ),
\end{eqnarray*}
where $T_1, \ldots, T_J$ denote the successive arrival times on $[0,K]$ of a standard Poisson process.  Thus,
$$
\lim_{\delta \rightarrow 0}P\left( \inf_{i=1, \ldots, J}\{T_i^{1/\gamma}-T_{i-1}^{1/\gamma}\} < \delta\right) = 0.
$$
and therefore \eqref{line:tight_condB} holds.  This completes the proof of (i).

Now we prove (ii):
Fix $0<c<\infty$.  We first write
\begin{eqnarray}
\sup_{0< x \leq ca_n} \left|\frac{\widehat f_n(x)}{f_0(x)}-1\right|
&=& \sup_{0< t \leq c} \left|\frac{na_n \widehat{f}_n(t a_n ) }{na_n f_0(t a_n )}-1\right| .
           \label{line:A1:calc1}
\end{eqnarray}
Suppose we could show that the ratio process
$na_n \widehat f_n(a_n t)/na_n f_0(a_n t )$ converges to the process
$t^{1-\gamma} \widehat h_\gamma(t)/\gamma$ in $D[0,\infty)$.
Then the conclusion follows by noting that the functional
$h\mapsto \sup_{0<t\leq c}|h|$ is continuous in the Skorohod topology
as long as $c$ is not a point of discontinuity of $h$ (\citet[Proposition VI 2.4, page 339]{MR1943877}).  Since $\NN(t^\gamma)$ is stochastically continuous
(i.e. $P( \NN(t^{\gamma} ) - \NN(t^{\gamma}-) > 0) = 0$   for each fixed $t>0$),
$t^{1-\gamma} \widehat h_\gamma(t)/\gamma$ is almost surely continuous at $c$.

It remains to prove convergence of the ratio.  Fix $K>c$, and again we may assume that
$K$ is a continuity point.  Consider first the term in the denominator, $na_n f_0(a_n t)$:
it follows from  (\ref{ConseqGnedenkoHypothesisForDensity}) that
\begin{eqnarray*}
g_n (t) \equiv (na_n f_0( a_n t))^{-1}
 \rightarrow  \gamma^{-1} t^{1-\gamma} \equiv g(t)
\end{eqnarray*}
where $g$ is monotone increasing and uniformly continuous on $[0,K]$.  Thus $g_n \rightarrow g$
in $C[0,K]$.  Since the term in the numerator satisfies
$h_n (t) \equiv na_n \widehat{f}_n (a_n t) \Rightarrow \widehat{h}_{\gamma}(t) \equiv h(t) $ in $D[0,K]$,
it follows that $g_n h_n \Rightarrow g h$ in $D[0,K]$, as required.
 Here, we have again used the continuity of the supremum.  This completes the proof of~(ii).
\hfill $\Box$
\medskip

Before proving Corollaries~\ref{cor:BoundedCase} - \ref{cor:PolyGrowthCase} we state the following lemma. 
\begin{lem}\label{lem:an}
Suppose that $a_n = p(1/n)$ for some function with $p(0) = 0$ satisfying 
$\lim_{x \rightarrow 0+} p'(x) f_0 (p(x)) = 1$.  Then $n F_0 (a_n ) \rightarrow 1$.  
\end{lem}
\par\noindent{\bf Proof:}
This follows easily from l'H\^{o}pital's rule, since 
\begin{eqnarray*}
\lim_{n\rightarrow \infty} n F_0(a_n) = \lim_{x\rightarrow 0+} \frac{F_0(p(x))}{x}= \lim_{x\rightarrow 0+} f_0(p(x))p'(x).
\end{eqnarray*}\hfill $\Box$

\par\noindent
{\bf Proof of Corollary~\ref{cor:BoundedCase}:}
Under the assumption {\bf G0} we see that $F_0 (x) \sim f_0 (0+) x$ as $x \rightarrow 0$,
so (\ref{GnedenkoHypothesis}) holds with $\gamma =1$. 
The claim that $a_n = 1/(n f_0 (0+))$ satisfies $n F_0 (a_n )\rightarrow 1$ follows from Lemma \ref{lem:an}
with $p(x) = x/f_0 (0+)$.
For (i) note that $\widehat{h}_1 (0) = \widehat{h}_1 (0+) = \sup_{t>0} ( \NN (t) /t)$,
and the indicated equality in distribution follows from \citet{MR0107315}; see
Proposition~\ref{DistributionOfYSubAlpha} and its proof.  (ii) follows directly from (i) of
Theorem~\ref{thm:GeneralTheoremGrenanderAtZero}.
To prove (iii), note that from (ii) of Theorem~\ref{thm:GeneralTheoremGrenanderAtZero}
it suffices to show that
\begin{eqnarray}
 \sup_{0<t\le c}\left| \widehat{h}_{1} (t) - 1\right | = \left|\widehat{h}_1 (0+)-1\right| = \widehat{h}_1 (0+)-1 = Y_1 -1
 \label{EqualityToProvePartZero}
 \end{eqnarray}
for each $c>0$
where $\widehat{h}_1 (t)$ is the right derivative of the LCM of $\NN(t)$.
The equality in (\ref{EqualityToProvePartZero}) holds if
$\widehat{h}_1 (c) >1$, since $\widehat{h}_1$ is decreasing by definition.
By the switching relation (\ref{ClosedSwitchingRelationGrenander}), we have the equivalence
\begin{eqnarray*}
\{\widehat{h}_1 (c)>1\} = \{\widehat{s}^L (1) > c\} .
\end{eqnarray*}
The  equality in (\ref{EqualityToProvePartZero}) thus follows if
$\widehat{s}^L(1)=\infty$.    That is, if
\begin{eqnarray*}
\NN(t)-t < \sup_{y\geq 0} \{\NN(y)-y\} \ \  \mbox{ for all finite }t.
\end{eqnarray*}
Let $W = \sup_{y\geq 0} \{\NN(y)-y\}$.
\citet[pages 570-571]{MR0107315} showed that $P(W \leq x)=0$ for $x \ge 0$; i.e. $P(W = \infty) = 1$.
\hfill $\Box$
\medskip

\par\noindent
{\bf Proof of Corollary~\ref{cor:LogGrowthCase}:}
Under the assumption {\bf G1} we see that $F_0 (x) \sim C_1 x (\log (1/x))^{\beta}$ as $x \rightarrow 0$,
so (\ref{GnedenkoHypothesis}) holds with $\gamma =1$. 
The claim that $a_n = 1/(C_1 n (\log n)^{\beta} )$ satisfies $n F_0 (a_n )\rightarrow 1$ follows from Lemma \ref{lem:an} with 
$p(x) = x/(C_1 \log (1/x))^{\beta}$.
For (i) note that $\widehat{h}_1 (0) = \widehat{h}_1 (0+) = \sup_{t>0} ( \NN (t) /t)$ just as in the proof of
Corollary~\ref{cor:BoundedCase}.
 (ii) again follows directly from (i) of
Theorem~\ref{thm:GeneralTheoremGrenanderAtZero}, and the proof of (iii) is just the same
as in the proof of Corollary~\ref{cor:BoundedCase}.
\hfill $\Box$
\medskip

\par\noindent
{\bf Proof of Corollary~\ref{cor:PolyGrowthCase}:}
Under the assumption {\bf G2} we see that $F_0 (x) \sim C_2 x^{1-\alpha} /(1-\alpha)$ as $x \rightarrow 0$,
so (\ref{GnedenkoHypothesis}) holds with $\gamma =1-\alpha$.  
The claim that $a_n =  \{ (1-\alpha) /(nC_2 )\}^{1/(1-\alpha)}$ satisfies $n F_0 (a_n ) \rightarrow 1$ follows from 
Lemma \ref{lem:an} with $p(x) = ( (1-\alpha)x/C_2)^{1/(1-\alpha)}$.
For (i) note that
$$
\widehat{h}_{1-\alpha} (0) = \widehat{h}_{1-\alpha} (0+)
= \sup_{t>0} ( \NN (t^{1-\alpha}) /t) = \sup_{s>0} ( \NN (s) / s^{1/(1-\alpha)} )
$$
much as in the proof of
Corollary~\ref{cor:BoundedCase}.
 (ii) and (iii) follow directly from (i) and (ii) of
Theorem~\ref{thm:GeneralTheoremGrenanderAtZero}.
\hfill $\Box$

\medskip

\par\noindent
{\bf Proof of Proposition~\ref{DistributionOfYSubAlpha}:}
        The part of the proposition with $\gamma =1$ follows from
        \citet[pages 570-571]{MR0107315}; this is closely related to a classical result of
        \citet{MR0012388} for the empirical distribution function; see e.g.
        \citet[Theorem 9.1.2, page 345]{MR838963}.

        The proof for the case $\gamma <1$ proceeds much along the lines of
        \citet[pages 103--105]{MR694539}.  Fix $x>0$ and $\gamma <1$. We aim
        at establishing an expression for the distribution function of
        $Y_{\gamma} \equiv \sup_{s>0}( \NN(s)/s^{1/\gamma} )$
        at $x>0$.   First, observe that
        \begin{eqnarray}
            P(Y_{\gamma} \leq x) 
            &= & P\left( \sup_{s>0}\left\{\frac{\NN(s)}{s^{1/\gamma}}\right\} \leq x\right)   \nonumber \\
            &= & P( \NN(t) \leq U(t)\ \quad\text{for all $t>0$})    \label{FirstEqualDFofYGamma}
       \end{eqnarray}
        where the 
        function $U(t) = xt^{1/\gamma}$. 
        For $j\in \NN$ let $t_j:=(j/x)^{\gamma}$, and note that 
        $t_1<t_2<~\ldots$ and $U(t_j) = j$.

        Define sets $B$ and $C$ by
        $$
        B \equiv [\NN(t_k)\neq k\,;
        \text{ for all $k\geq 1$}] \ \ \mbox{and} \ \
        C \equiv [\NN(s) > U(s)\,;
        \text{ for some $s>0$}].
        $$
        Then
         $P(B\cap C)=0$
        as a consequence of the following argument:
        Suppose that 
        there exists some $t>0$ and $k\in\NN$ such that
        $k=\NN(t)>U(t)$ and $\NN(t_i)\neq i$, for all $i\geq 1$.
        It then follows that $t_k>t$, for otherwise it follows
        that $k=U(t_k)\leq U(t)$, as $U(\cdot)$ is
        increasing, which is a contradiction. Therefore, $t_k >t$
        implies that $\NN(t_k)>\NN(t)=k$, as $\NN(\cdot)$ is
        non--decreasing while $\NN(t_k)=k$ is disallowed, by
        hypothesis. Hence, $\NN(t_i) > i$ holds true for all $i\geq k$,
        for otherwise there would exist some $j\geq k$ such that
        $\NN(t_j)=j$, since $\NN(\cdot)$ is a counting process.
        Therefore, for each $i\geq k$ we have that $\NN(s)\geq i+1$
        holds for all $t_i\leq s\leq t_{i+1}$ and, consequently, that
        $\NN(s) \geq U(s)$ holds for all $s\geq t_k$. This implies
        that $B\cap C \subseteq [ \liminf_{s\to\infty}\{ \NN(s) /
        s^{1/\gamma}\}\geq x]$ and therefore $P(B\cap C) = 0$, since the SLLN
        implies that $\NN(s)/s^{1/\gamma} \rightarrow 0$ holds almost
        surely, for fixed $\gamma < 1$. We thus conclude that
        $P(B\cap C) = 0$.

      We conclude that $P(C) = P(C \cap B^c)$.  Furthermore, since $U$ is a strictly increasing
      function, and since $\NN$ has jumps at the points $\{t_k\}$ with probability zero, 
      we also find that $P(C \cap B^c) = P(B^c)$.  
        Finally, partition $B^c$ as $B^c = \cup_{k=1}^{\infty} A_k$
        for the disjoint sets 
        $A_k\equiv [\NN(t_k) = k, \NN(t_j)\neq j \text{ for all $1\leq j < k$}]$, 
        $k\geq 1$. Combining all
        arguments above, we conclude that
        $$
        P(Y_{\gamma} \leq x) = 1-P(C) = 1-\sum_{k=1}^{\infty} P(A_k)
        $$
        where $P(A_1 ) = P( \NN(t_1) = 1) = p(t_1; 1)$, and, for $k\ge 2$,
        $P(A_k)$ may be written as
        \begin{eqnarray*}
        \lefteqn{P(\NN(t_k) = k) - P( \{ \NN(t_k) = k \} \cap \{ \NN(t_i) \not= i, \, i < k \}^c )}\\
        & = & P( \NN(t_k) = k) - \sum_{j=1}^{k-1} P(\NN(t_k ) = k, \, \NN (t_j ) = j, \, \NN(t_i) \not= i, \, i < j) \\
        & = & P(\NN(t_k ) = k ) - \sum_{j=1}^{k-1} P( \NN(t_k) - \NN(t_j) = k-j) P( \NN( t_j ) = j, \NN(t_i) \not= i, \, i< j) .
        \end{eqnarray*}
        The result follows.
        \hfill $\Box$


\setcounter{section}{6}
\section*{Appendix B: Definitions from Convex Analysis}

The \textit{epigraph} (\textit{hypograph}) of a function $f$ from a subset $S$ of $\RR^d$ to $[-\infty, +\infty]$ is the subset 
$\mbox{epi}(f)$ ($\mbox{hypo}(f)$) of $\RR^{d+1}$ defined by
\begin{eqnarray*}
\mbox{epi}(f) &= &\{ (x,t): \, x \in S, \, t \in \RR, \, t \ge f(x) \},\\
\mbox{hypo}(f)&=& \{ (x,t) : \ x \in S, \, t \in \RR; \ t \le f(x) \}.
\end{eqnarray*}
The function $f$ is \textit{ convex} if $\mbox{epi}(f)$ is a convex set.
The \textit{ effective domain} of a convex function $f$ on $S$ is 
\begin{eqnarray*}
\mbox{dom} (f) = \{ x \in \RR^d : \, (x,t) \in \mbox{epi}(f) \ \mbox{for some } \ t \} = \{ x \in \RR^d : \ f(x) < \infty \} .
\end{eqnarray*}

The $t-$sublevel set of a convex function $f$ is the set
$
C_t = \{ x \in \mbox{dom}(f) : \ f(x) \le t \}, 
$
and the $t-$superlevel set of a concave function $g$ is the set 
$
S_t = \{ x \in \mbox{dom}(g) : \ g(x) \ge t \} .
$
The sets $C_t$, $S_t$ are convex.
The \textit{ convex hull} of a set $S \subset \RR^d$, denoted by $\mbox{conv}(S)$,
 is the intersection of all the convex sets containing $S$.

A convex function $f$ is said to be \textit{ proper} if its epigraph is non-empty and contains no vertical lines; i.e. if $f(x) < +\infty$ 
for at least one $x$ and $f(x) > - \infty$ for every $x$.  
Similarly, a concave function $g$ is \textit{ proper} if the convex function $-g$ is proper.
The \textit{ closure of a concave function} $g$, denoted by $\mbox{cl}(g)$, is the pointwise infimum of all affine functions
$h\ge g$.  If $g$ is proper, then
$$
\mbox{cl} (g) (x) = \limsup_{y \rightarrow x} g(y) .
$$
For every proper convex function $f$ there exists closed
proper convex function $\cl(f)$ such that $\mbox{epi} (\cl(f))\equiv\cl(\mbox{epi} (f))$. 
The \textit{ conjugate function}  $g^*$ of a concave function $g$ is defined by
$$
g^* (y) = \inf \{ \langle x, y \rangle - g(x) : \ x \in \RR^d \} ,
$$
and the conjugate function $f^*$ of a convex function $f$ is defined by 
$$
f^* (y) = \sup \{ \langle x,y \rangle - f(x) : \ x \in \RR^d \} .
$$
If $g$ is concave, then $f = -g$ is convex and  $f$ has conjugate $f^* (y) = - g^* (-y)$.

A \textit{ complete non-decreasing curve} is a subset of $\RR^2$ of the form
$$
\Gamma = \{ (x,y) : \ x \in \RR,  \ y \in \RR, \ \varphi_{-} (x) \le y \le \varphi_+ (x) \} 
$$
for some non-decreasing function $\varphi$ from $\RR$ to $[-\infty, +\infty]$ which is not everywhere infinite.
Here $\varphi_{+}$ and $\varphi_{-}$ denote the right and left continuous versions of $\varphi$
respectively.  A vector $y \in \RR^d$ is said to be a \textit{ subgradient} of a convex function $f$ at a point
$x$ if 
$$
f(z) \ge f(x) + \langle y, z-z \rangle \qquad \mbox{for all} \ \ z \in \RR^d .
$$
The set of all subgradients of $f$ at $x$ is called the \textit{ subdifferential of $f$ at $x$}, 
and is denoted by $\partial f (x)$.

A \textit{ face} of a convex set $C$ is a convex subset $B$ of $C$ such that every closed line segment in $C$ 
with a relative interior point in $B$ has both endpoints in $B$.  
If $B$ is the set of points where a linear function $h$ achieves its maximum over $C$, then $B$ is a face of $C$. 
If the maximum is achieved on the relative interior of a line segment $L \subset C$, then $h$ must be constant on $L$
and $L \subset B$.  A face $B$ of this type is called an \textit{ exposed face}.


\vskip .65cm
\noindent
Centre de Recherche en Math\'{e}matiques de la D\'{e}cision\\
Universit\'{e} Paris-Dauphine,  Paris, France
\vskip 2pt
\noindent
E-mail: fadoua@ceremade.dauphine.fr
\vskip .3cm
\noindent
Department of Mathematics and Statistics\\
York University, Toronto, Canada
\vskip 2pt
\noindent
E-mail: hkj@mathstat.yorku.ca
\vskip .3cm
\noindent
Department of Mechanical Engineering\\
Frederick University Cyprus,  Nicosia, Cyprus
\vskip 2pt
\noindent
E-mail: m.pavlides@frederick.ac.cy
\vskip .3cm
\noindent
Department of Statistics\\
University of Washington,  Seattle, USA
\vskip 2pt
\noindent
E-mail: arseni@stat.washington.edu
\vskip .3cm 
\noindent
Department of Statistics\\
University of Washington,  Seattle, USA
\vskip 2pt
\noindent
E-mail: jaw@stat.washington.edu\\

\end{document}